\documentclass[final]{siamltex}

\usepackage{graphicx}
\usepackage{epstopdf}

\usepackage{amsmath}
\usepackage{amssymb}

\usepackage{tikz}
\usepgflibrary{plotmarks}
\usetikzlibrary{plotmarks}
\usetikzlibrary{shapes,shadows}

\usepackage{bm}
\usepackage{subfigure}

\usepackage{url}
\usepackage{xcolor}
\definecolor{lightgray}{gray}{0.5}
\definecolor{darkred}{rgb}{.6,0,0}
\definecolor{darkblue}{rgb}{0,0,0.5}
\usepackage[colorlinks=true,citecolor=darkblue,urlcolor=darkblue,linkcolor=darkred]{hyperref}

\renewcommand{\thefigure}{\arabic{section}.\arabic{figure}}
\makeatletter \@addtoreset{figure}{section} \makeatletter

\renewcommand{\thetable}{\arabic{section}.\arabic{table}}
\makeatletter \@addtoreset{table}{section} \makeatletter

\newtheorem{thm}{Theorem}[section]
\newtheorem{prop}[thm]{Proposition}

\newtheorem{alg}[thm]{Algorithm}

\newtheorem{remark}[thm]{Remark}

\numberwithin{equation}{section}

\begin{document}


\title{The Closest Point Method and multigrid solvers for \\ elliptic equations on surfaces}
\author{Yujia Chen \and Colin B.~Macdonald\thanks{Mathematical Institute, University of Oxford, OX1 3LB, UK (\url{cheny1,macdonald@maths.ox.ac.uk}). This publication is based on work supported by Award No. KUK-C1-013-04, made by King Abdullah University of Science and Technology (KAUST).}}

\maketitle







\begin{abstract}
Elliptic partial differential equations are important both from application and analysis points of views.
In this paper we apply the Closest Point Method to solve elliptic equations on general curved surfaces.
Based on the closest point representation of the underlying surface, we formulate an
embedding equation for the surface elliptic problem, then discretize it
using standard finite differences and interpolation schemes
on banded, but uniform Cartesian grids. We prove the convergence of the
difference scheme for the Poisson's equation on a smooth closed curve. In order to solve the resulting
large sparse linear systems, we propose a specific geometric multigrid method in the setting
of the Closest Point Method. Convergence studies both in the accuracy of the
difference scheme and the speed of the multigrid algorithm show that our
approaches are effective.
\end{abstract}

\begin{keywords}
  Surface elliptic problem, Surface Poisson problem, Laplace--Beltrami
  operator, Closest Point Method, Geometric multigrid method.
\end{keywords}

\begin{AMS}
58J05, 65N06, 65N55
\end{AMS}

\pagestyle{myheadings}
\thispagestyle{plain}
\markboth{Y.-J. CHEN AND C. B. MACDONALD}{CLOSEST POINT METHOD FOR ELLIPTIC SURFACE PROBLEMS}








\section{Introduction}

Partial differential equations (PDEs) on surfaces are important for many applications
in engineering \cite{plain:surfactant}, material science \cite{plain:phase-change},
biology \cite{plain:wound-healing}, and image processing \cite{plain:segmentation}.
In all these applications, numerical algorithms are essential for solving PDEs.
There are a variety of techniques for solving surface PDEs numerically, typically depending on the geometrical representation of the underlying surface:
\begin{itemize}
\item One can find a parametrization of the surface \cite{plain:param}, formulate and solve the PDE in the parameter space.
\item One can also find a triangulation or a polygonal approximation of the surface, and solve the PDE on the
approximated surface by finite differences \cite{plain:texture} or finite elements \cite{plain:sfem-error,plain:sfem}.
\item Based on a level set representation of the surface, one can formulate an embedding PDE corresponding to the surface PDE,
then use either finite differences \cite{plain:Bertalmio,plain:Greer} or finite elements \cite{plain:Burger,plain:dc1,plain:Deckelnick}
to numerically solve the PDE in the embedding space. This level set approach is desirable for PDEs on moving surfaces
\cite{plain:Zhao,plain:dc2}.
\item Based on a closest point representation of the surface, one can also formulate an embedding equation (in a different way from the
level set approach), then solve the embedding equation in a banded, uniform Cartesian grid with finite differences
\cite{plain:RuuthMerriman,plain:lscpm,plain:icpm}. This approach is called the Closest Point Method.
\item Starting from a point cloud (supposed to be sampled from the surface), one can reconstruct the surface by radial basis functions,
and solve the PDE on the reconstructed surface via the radial basis function finite difference method \cite{plain:OGr,plain:kernel}.
There are also other approaches to compute on point clouds with no
reconstruction procedure \cite{MacdonaldMerrimanRuuth:ptclouds} or with only
local reconstruction \cite{LiangZhao:ptclouds}.
\end{itemize}

In this paper we will adopt the Closest Point Method, and focus on solving elliptic PDEs on surfaces.
Similarly to the level set approach, the Closest Point Method extends the surface PDE to the embedding space of the surface
in such a way that the embedding equation agrees with the original PDE on the surface. However, the Closest Point Method
extends the solution via the \emph{closest point extension} so that the extended solution is constant along the normals to the surface,
and thus the Cartesian differential operators agree with the surface intrinsic differential operators on the surface.
One does not need to modify the Cartesian differential operators as in the case of the level set approach.

The Closest Point Method has been successfully applied to solving a wide variety of time-dependent PDEs
\cite{plain:RuuthMerriman,plain:lscpm,plain:icpm} and eigenvalue problems for the Laplace--Beltrami operator \cite{plain:eigen}
on surfaces, but less attention has been paid to solving elliptic problems.
Motivated by constructing implicit schemes for the in-surface heat equation,
Macdonald and Ruuth \cite{plain:icpm} propose a way to discretize the embedding operator
corresponding to the Laplace--Beltrami operator;
this idea is also successfully applied to solving eigenvalue problems for
the Laplace--Beltrami operator on surfaces \cite{plain:eigen}.
One can apply the approach in \cite{plain:icpm,plain:eigen} to deal with the elliptic operators,
but this requires a properly extended right-hand side function from the surface to the embedding space.
A simple closest point extension of the right-hand side seems to work in
practice but the resulting schemes are hard to analyze because the
solutions are \emph{not} constant in the direction of the normal of
the surface.
We discuss this further in Section~\ref{section-inconsistency}.

To overcome this issue, we formulate an embedding equation for the surface elliptic problem
using a different approach from \cite{plain:icpm,plain:eigen}.
The new embedding equation agrees with the original PDE on the surface and
has a solution which is constant along the normals to the surface.
Similar ideas have been investigated for method-of-lines approaches to solve time-dependent
problems involving more general surface differential operators \cite{plain:mol},
but here we focus on elliptic problems.
We then construct a finite difference scheme which is consistent with the embedding equation.
Furthermore, for the shifted Poisson's equation on a smooth closed curve embedded in $\mathbb{R}^2$,
we prove that the difference scheme
is second-order convergent in the $\infty$-norm under some smoothness assumptions.

The difference scheme yields a sparse linear system.
If $\mathcal{S}$ is a curve in $\mathbb{R}^2$, the resulting linear system can be easily solved
by sparse direct solvers such as Matlab's backslash; but if $\mathcal{S}$ is a surface in $\mathbb{R}^3$
or higher-dimensional spaces,
efficient iterative solvers are needed.
We propose a specific geometric multigrid algorithm in the setting of the Closest Point Method,
making full use of the closest point representation of the surface
and uniform Cartesian grids surrounding the surface.

The rest of the paper is organized as follows. In Section~\ref{section-surface-poisson} we
formulate the embedding equation for the surface elliptic problem. In Section~\ref{section-discretization} we discretize
the embedding equation for the surface Poisson problem, and analyze the difference scheme;
we also provide some numerical examples for Poisson equations on closed curves in $\mathbb{R}^2$ to validate our analysis.
In Section~\ref{section-v-cycle}, motivated by solving the large sparse linear systems arising from
problems where the embedding space is $\mathbb{R}^3$ or higher, we propose a V-Cycle Multigrid algorithm
based on the closest point representation of the surface.
We then present examples on various surfaces and with variable diffusion coefficients to show the efficiency of our specific V-Cycle Multigrid algorithm.
In Section~\ref{section-conclusion} we draw several conclusions and discuss some future work.

\section{Continuous embedding equation}
\label{section-surface-poisson}
\subsection{Embedding equation for surface Poisson problem}
\label{section-poisson}
For ease of derivation, we first focus on the shifted Poisson's equation
on a smooth closed surface $\mathcal{S}$:
\begin{equation}
\label{goal}
-\Delta_{\mathcal{S}}u(\bm{y})+cu(\bm{y}) = f(\bm{y}), \qquad \bm{y}\in \mathcal{S},
\end{equation}
where $f(\bm{y})$ is a continuous function defined on the surface $\mathcal{S}$,
and we assume $c$ to be a positive constant for simplicity.\footnote{
Our derivations and proofs also work for the case where $c(\bm{y})$ is a surface function satisfying
 $c(\bm{y})\geq C>0$, where $C$ is a positive constant.}
Since $c>0$, (\ref{goal}) has a unique solution.

Because $\mathcal{S}$ is smooth, there is a tubular neighborhood $B(\mathcal{S})$ surrounding
$\mathcal{S}$ \cite{plain:hirsch,plain:cpfunctions} such that
the following function
is well defined in $B(\mathcal{S})$.
We define the (Euclidean) closest point function
$\text{cp}:B(\mathcal{S})\rightarrow\mathcal{S}$ such that
for each point $\bm{x}\in B(\mathcal{S})$,
$\text{cp}(\bm{x})$ is a surface point which is closest to
$\bm{x}$ in Euclidean distance \cite{plain:RuuthMerriman}.
We want to define an embedding equation in $B(\mathcal{S})$, whose solution agrees with
the solution of the surface PDE \eqref{goal} on $\mathcal{S}$.

We construct the right-hand side function $\tilde{f}$ in $B(\mathcal{S})$ via the closest point extension,
\begin{equation*}
\tilde{f}(\bm{x}) = f(\text{cp}(\bm{x})), \qquad \bm{x}\in B(\mathcal{S}).
\end{equation*}
From this construction of $\tilde{f}$  and
the idempotence of the closest point extension, we have that
\begin{equation}
\label{normal-rhs-con}
\tilde{f}(\bm{x}) = \tilde{f}(\text{cp}(\bm{x})), \qquad \bm{x}\in B(\mathcal{S}).
\end{equation}
We also want to find a solution $\tilde{u}$ in the band $B(\mathcal{S})$ which is constant along the normal direction to the surface,
i.e., we will need to impose this condition on $\tilde{u}$,
\begin{equation}
\label{normal-u-con}
\tilde{u}(\bm{x}) = \tilde{u}(\text{cp}(\bm{x})),  \qquad \bm{x}\in B(\mathcal{S}).
\end{equation}

Since $\tilde{u}$ is constant along the normal direction to the surface,
the Cartesian Laplacian of $\tilde{u}$ equals the Laplace--Beltrami function of $\tilde{u}$ on the surface \cite{plain:RuuthMerriman},
and thus
\begin{displaymath}
-\Delta \tilde{u}(\bm{y}) + c\tilde{u}(\bm{y}) = \tilde{f}(\bm{y}), \qquad \bm{y}\in\mathcal{S}.
\end{displaymath}
Replacing $\bm{y}$ with $\text{cp}(\bm{x})$, 
we have
\begin{displaymath}
-[\Delta \tilde{u}](\text{cp}(\bm{x})) + c\tilde{u}(\text{cp}(\bm{x})) = \tilde{f}(\text{cp}(\bm{x})), \qquad \bm{x}\in\text{B}(\mathcal{S}).
\end{displaymath}
Here the square bracket $[\cdot]$ surrounding $\Delta\tilde{u}$ means that we first calculate the Cartesian Laplacian of $\tilde{u}(\bm{x})$,
and then evaluate at the closest point of $\bm{x}$ on the surface.
By assumption, both $\tilde{f}$ and $\tilde{u}$ are closest point extensions
((\ref{normal-rhs-con}) and (\ref{normal-u-con})), so we have
\begin{subequations}
\label{pde-constraint-system}
\begin{align}
\label{pde-band}
-[\Delta\tilde{u}](\text{cp}(\bm{x})) + c\tilde{u}(\bm{x}) = \tilde{f}(\bm{x}), \\
\label{constraint}
\text{subject to}\quad \tilde{u}(\bm{x}) = \tilde{u}(\text{cp}(\bm{x})), \quad
\bm{x}\in B(\mathcal{S}).
\end{align}
\end{subequations}

Equation (\ref{pde-band}) ensures that the embedding equation agrees with the original surface PDE on the surface,
and we will call it the \emph{consistency condition}.
Equation (\ref{constraint}) forces $\tilde{u}$ to be constant along the normal direction to the surface
(so that we can replace the surface Laplacian with the Cartesian Laplacian), and we will call it the \emph{side condition}.
Combining the two equations (\ref{pde-band}) and (\ref{constraint}), we get the embedding equation,
\begin{equation}
\label{embedding}
-[\Delta\tilde{u}](\text{cp}(\bm{x})) + c\tilde{u}(\bm{x}) + \gamma(\tilde{u}(\bm{x})-\tilde{u}(\text{cp}(\bm{x})))
= \tilde{f}(\bm{x}), \quad \bm{x}\in B(\mathcal{S}),
\end{equation}
where $\gamma$ is a parameter. As we shall show next, the choice of $\gamma$ does not have much impact on the continuous equation (\ref{embedding}); but later in our numerical scheme, $\gamma$ plays a role in balancing the side and consistency conditions.

So far we have shown that the system (\ref{pde-constraint-system}) leads to equation (\ref{embedding}).
Following the ideas of \cite{plain:mol},
 we are going to show that in the continuous setting and if $\gamma\not=-c$, equation (\ref{embedding}) implies \emph{both} equations (\ref{pde-band}) and (\ref{constraint}).
We rearrange equation (\ref{embedding}) as
\begin{equation}
\label{rearrange-embedding}
-[\Delta\tilde{u}+\gamma\tilde{u}](\text{cp}(\bm{x})) = \tilde{f}(\bm{x}) - (\gamma+c) \tilde{u}(\bm{x}), \quad \bm{x}\in B(\mathcal{S}).
\end{equation}
Here the square bracket $[\cdot]$ surrounding $\Delta\tilde{u}+\gamma\tilde{u}$
means that we first calculate the Cartesian Laplacian of $\tilde{u}(\bm{x})$
and add the result with $\gamma\tilde{u}(\bm{x})$, then evaluate the function
$\Delta\tilde{u}+\gamma\tilde{u}$ at the closest point of $\bm{x}$ on the surface.
The left-hand side of equation (\ref{rearrange-embedding}) is constant along the direction normal to $\mathcal{S}$,
so the right-hand side $\tilde{f}(\bm{x}) - (\gamma+c)\tilde{u}(\bm{x})$ also has to be constant along the normal direction to $\mathcal{S}$,
i.e.,
\begin{displaymath}
\tilde{f}(\bm{x}) - (\gamma+c)\tilde{u}(\bm{x})=\tilde{f}(\text{cp}(\bm{x})) - (\gamma+c)\tilde{u}(\text{cp}(\bm{x})).
\end{displaymath}
Since $\tilde{f}(\bm{x})=\tilde{f}(\text{cp}(\bm{x}))$, if $\gamma\not=-c$, we have that
 $\tilde{u}(\bm{x})=\tilde{u}(\text{cp}(\bm{x}))$, i.e., equation (\ref{constraint}) holds.
Substituting (\ref{constraint}) back into equation (\ref{embedding}), we have that equation (\ref{pde-band}) holds.
In other words, from equation (\ref{embedding}) with $\gamma\not=-c$,
we derive that both equations (\ref{pde-band}) and (\ref{constraint}) hold.

In summary, our logic is as follows: in order to solve the original surface PDE (\ref{goal}),
we seek $\tilde{u}(\bm{x})$ in $B(\mathcal{S})$ satisfying both the consistency condition (\ref{pde-band}) and
the side condition (\ref{constraint}), which is in turn equivalent to equation (\ref{embedding}) (provided $\gamma\not=-c$).
So from now on we can focus on equation (\ref{embedding}).

\subsection{More general elliptic operators}
\label{sec:general_ops}
The above ideas extend naturally to more general linear surface elliptic PDEs,
\begin{equation}
\label{original-general-A}
-\nabla_{\mathcal{S}}\cdot(\bm{A}(\bm{y})\nabla_{\mathcal{S}} u(\bm{y}))
+ c(\bm{y})u(\bm{y}) = f(\bm{y}), \quad \bm{y}\in\mathcal{S},
\end{equation}
where $\bm{A}(\bm{y})$ is a symmetric and semi-positive definite matrix.

Again we construct $\tilde{f}(\bm{x})=f(\text{cp}(\bm{x}))$, $\bm{x}\in B(\mathcal{S})$, and would like to enforce the side condition
$\tilde{u}(\bm{x}) = \tilde{u}(\text{cp}(\bm{x}))$ for the solution $\tilde{u}$ defined in $B(\mathcal{S})$. Using the generalized closest point principle (\cite[Theorem 4.3]{plain:cpfunctions}),
the following equations agree with surface PDE (\ref{original-general-A}) on the surface,
\begin{subequations}\label{pde-constraint-general}
\begin{align}
\label{pde-band-general}
-[\nabla\cdot(\bm{A}(\text{cp}(\bm{x}))\nabla\tilde{u})](\text{cp}(\bm{x})) + c\tilde{u}(\bm{x}) = \tilde{f}(\bm{x}), \\ 
\label{constraint-general}
\text{subject to}\quad \tilde{u}(\bm{x}) = \tilde{u}(\text{cp}(\bm{x})), \quad
\bm{x}\in B(\mathcal{S}).
\end{align}
\end{subequations}
Here the square bracket $[\cdot]$ surrounding $\nabla\cdot(\bm{A}(\text{cp}(\bm{x}))\nabla\tilde{u})$ means that we first perform the gradient on $\tilde{u}$,
multiply it by $\bm{A}(\text{cp}(\bm{x}))$, then take the divergence, and finally evaluate the results at $\text{cp}(\bm{x})$.
Adding (\ref{pde-band-general}) and (\ref{constraint-general}) together, we have
\begin{equation}
\label{embedding-general}
-[\nabla\cdot(\bm{A}(\text{cp}(\bm{x}))\nabla\tilde{u})](\text{cp}(\bm{x})) + c\tilde{u}(\bm{x}) + \gamma(\tilde{u}(\bm{x})-\tilde{u}(\text{cp}(\bm{x})))
= \tilde{f}(\bm{x}), \quad \bm{x}\in B(\mathcal{S}),
\end{equation}

Similarly to the Poisson case in Section~\ref{section-poisson}, we have the following equivalence relationship.
\begin{thm}
\begin{enumerate}
\item If $\tilde{u}\in B(\mathcal{S})$ is a solution of (\ref{pde-constraint-general}), then
$\tilde{u}$ is a solution of (\ref{embedding-general}).
\item If $\tilde{u}\in B(\mathcal{S})$ is a solution of (\ref{embedding-general}), and $\gamma\not=-c$, then $\tilde{u}$
is a solution of (\ref{pde-constraint-general}).
\end{enumerate}
\end{thm}
{\it Proof.} The ideas are similar to the Poisson case in Section~\ref{section-poisson}, the only difference is that we replace
$\Delta\tilde{u}$ with $\nabla\cdot(\bm{A}(\text{cp}(\bm{x}))\nabla\tilde{u})$.

\bigskip
We focus our derivation and analysis on the surface Poisson problem via
the embedding equation (\ref{embedding}).
In most cases it is straightforward to extend our approach to more
general elliptic operators (e.g., we look at variable diffusion
coefficients in Section~\ref{section-sphere-vc}).

As we shall see later in our numerical discretization, by choosing suitable $\gamma$ in (\ref{embedding}) we can get a consistent and stable numerical scheme.
From (\ref{embedding}), no boundary conditions for the Laplacian operator are needed at the boundary of the band $B(\mathcal{S})$,
since we only need to evaluate $\Delta\tilde{u}$ at the closest points on the surface. We simply need a band large enough to enable this evaluation.

\subsection{A comparison with the Macdonald--Brandman--Ruuth approach}
\label{section-inconsistency}


If we apply the approach in \cite{plain:icpm,plain:eigen} to the
Laplace--Beltrami operator in (\ref{goal}), and
extend the right hand side via the closest point extension, then we get the following continuous embedding equation,
\begin{equation}
\label{embedding-icpm}
-\Delta(\tilde{u}(\text{cp}(\bm{x}))) + c\tilde{u}(\bm{x}) + \gamma(\tilde{u}(\bm{x})-\tilde{u}(\text{cp}(\bm{x}))) = \tilde{f}(\bm{x}),\quad \bm{x}\in B(\mathcal{S}).
\end{equation}
Note that the only difference between (\ref{embedding}) and (\ref{embedding-icpm}) is that in (\ref{embedding}) we first perform the Cartesian Laplacian operator on
$\tilde{u}$ and then evaluate at $\text{cp}(\bm{x})$, while in (\ref{embedding-icpm}) we first do the closest point extension of $\tilde{u}$ and then perform the Laplacian.
One can show (similarly to \cite{plain:eigen}), that equation
(\ref{embedding-icpm}) agrees with the original surface PDE
(\ref{goal}) on the surface,
and it has a unique solution provided that $\gamma\not=-c$.
We can re-arrange (\ref{embedding-icpm}), solving for $\tilde{u}(\bm{x})$, to obtain
\begin{equation*}
\tilde{u}(\bm{x}) = \frac{1}{c+\gamma}\big( \tilde{f}(\bm{x}) + \gamma\tilde{u}(\text{cp}(\bm{x})) + \Delta(\tilde{u}(\text{cp}(\bm{x}))) \big).
\end{equation*}
From this, $\tilde{u}$ is \emph{not constant along the normals to
  $\mathcal{S}$} because $\Delta(\tilde{u}(\text{cp}(\bm{x})))$ is not
constant along the normals.
This makes analysis of the method difficult.
Furthermore, it would also be difficult to make $\tilde{u}$ constant along the normals to $\mathcal{S}$ by modifying the
right-hand side (e.g., via some other method of extension).

In contrast, our new embedding equation (\ref{embedding}) not only agrees with the original surface PDE (\ref{goal}) on the surface,
but also has a solution $\tilde{u}$ which is constant in the normal direction to $\mathcal{S}$. The latter property is crucial to the consistency and stability analysis
for our numerical discretization in the next section.
In addition to these advantages for analysis, (\ref{embedding}) can lead to a consistent discretization with a sparser coefficient matrix than the one resulting from (\ref{embedding-icpm});
we discuss this further in Remark~\ref{remark-trunc}.
Furthermore, as we shall see in Section~\ref{section-v-cycle}, the construction of our multigrid solver relies on the property of $\tilde{u}$ being constant along the normal direction
to $\mathcal{S}$; it would not be straightforward to design a multigrid solver starting from (\ref{embedding-icpm}) since it gives a solution that is not constant along the normals.
Finally, as noted in Section~\ref{sec:general_ops}, our new approach extends naturally to non-constant-coefficient problems.

%
%
%

\section{Numerical discretization}
\label{section-discretization}
%
Similar to \cite{plain:icpm}, we derive a matrix formulation of a finite difference scheme for the embedding equation (\ref{embedding}).

\subsection{Construction of the difference scheme}
\label{sec:diffscheme}
Suppose the surface $\mathcal{S}$ is embedded in $\mathbb{R}^d$, and a uniform Cartesian grid is placed in the embedding space.
When discretizing $u(\text{cp}(\bm{x}))$, we need to assign the value at each grid point $\bm{x}$ to be the value at the corresponding
closest point $\text{cp}(\bm{x})$; since the closest point is generally not a grid point, its value is obtained through interpolation
of values at the surrounding grid points.
Similarly to \cite{plain:icpm}, we use tensor product Barycentric Lagrange interpolation \cite{plain:bc}, where
the interpolated value is a linear combination of the values at neighboring grid points in a hypercube.
We call these grid points neighboring the closest point $\text{cp}(\bm{x})$ the {\it{interpolation stencil}} of the scheme.

We also use standard finite difference schemes (e.g., the $\frac{1}{\Delta x^2}\{-4, 1, 1, 1, 1\}$ rule for the Laplacian in 2D),
and this requires a \emph{difference stencil} for each grid point.

\begin{figure}
\begin{center}
\begin{minipage}{0.5\textwidth}
\caption[Illustration of the computational band]{\label{double-band}
Examples of the lists of grid points $B_{interp}$ (indicated by {\tiny $\bullet$})
and $B_{edge}$ (indicated by {\tiny $\circ$}) where the curve $\mathcal{S}$ is
a circle. The interpolation stencil is a $4\times 4$ grid (arising, for
example, from using degree $3$ barycentric Lagrange interpolation) and
the shaded regions illustrate the use of this stencil at the three points
on the circle indicated by $\diamond$. Five-point difference stencils are
shown for two example points in $B_{interp}$, in one case illustrating the use
of points in $B_{edge}$. }
\end{minipage} \hspace{0.01\textwidth}
\begin{minipage}{0.32\textwidth}
\includegraphics[width=1\textwidth]{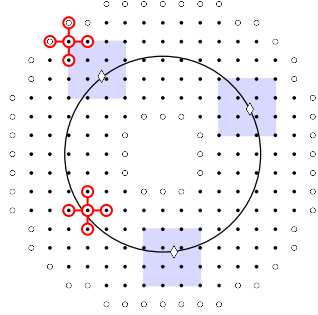}
\end{minipage}
\end{center}
\end{figure}

Let us review the definition of the two lists of discrete points in \cite{plain:icpm}.
The first list $B_{interp}=\{\bm{x}_1,\bm{x}_2,...,\bm{x}_m\}$ contains every grid point
which can appear in the interpolation stencil for some point on the surface $\mathcal{S}$.
The second one $B_{edge}=\{\bm{x}_{m+1},\bm{x}_{m+2},...,\bm{x}_{m+m_e}\}$
contains points which are not in $B_{interp}$ but appear in the difference stencil of some point in $B_{interp}$.
Figure~\ref{double-band}  shows an
example of the two lists $B_{interp}$ and $B_{edge}$ where the curve $\mathcal{S}$ is a
circle embedded in $\mathbb{R}^2$.
Let $B$ be a set of grid points such that $B_{interp} \cup B_{edge} \subseteq B$.
Let $n=\|B\|$ be the number of points in $B$.
We introduce the vector $u^h\in\mathbb{R}^{n}$ with entries $u_i\approx \tilde{u}(\bm{x}_i)$
for each grid point $\bm{x}_i\in B$,
and similarly, $f^h\in\mathbb{R}^{n}$ with entries $f_i=\tilde{f}(\bm{x}_i)$.

There are two closest point extensions in (\ref{embedding});
we can use different degrees of interpolation for the closest point evaluations of $\Delta\tilde{u}$ and $\tilde{u}$.
If we use degree $p$ interpolation for $[\Delta\tilde{u}](\text{cp}(\bm{x}))$,
and degree $q$ interpolation for $\tilde{u}(\text{cp}(\bm{x}))$, 
then for each grid point $\bm{x}_i\in B$, we have the following discretized equation,
\begin{equation}
\label{pointwise}
-\sum_{\bm{x}_j\in\text{I}_p(\text{cp}(\bm{x}_i))} \omega_j
\Bigg(\sum_{\bm{x}_k\in\text{Diff}(\bm{x}_j)} l_k u_k\Bigg) + cu_i
+ \gamma\Bigg(u_i- \sum_{\bm{x}_j\in\text{I}_q(\text{cp}(\bm{x}_i))}\bar{\omega}_j u_j\Bigg) = f_i.
\end{equation}
Here $\text{I}_p(\text{cp}(\bm{x}_i))$ is the degree $p$ interpolation stencil for point $\text{cp}(\bm{x}_i)$,
$\text{Diff}(\bm{x}_j)$ is the difference stencil for point $x_j$, and
$\text{I}_q(\text{cp}(\bm{x}_i))$ is the degree $q$ interpolation stencil for point $\text{cp}(\bm{x}_i)$;
$\omega_j$ are the weights of the interpolation scheme for evaluating $[\Delta\tilde{u}](\text{cp}(\bm{x}))$,
$\bar{\omega}_j$ are the weights of the interpolation scheme for evaluating $\tilde{u}(\text{cp}(\bm{x}))$,
and $l_k$ are the weights for the difference scheme approximating the Laplacian operator.

Since the interpolation and difference schemes are linear combinations of the entries of $u^h$,
we can write equation (\ref{pointwise}) in matrix form.
Denote the discrete solution by a column vector $u^h$, the discrete right-hand side function by a vector $f^h$,
the identity matrix by $I^h$,
the closest point interpolation matrices of degree $p$ and $q$ by $E_p^h$ and $E_q^h$,
and the discrete Laplace matrix by $L^h$.
The matrix form of the difference scheme for the embedding equation (\ref{embedding}) is
\begin{equation}
\label{discrete}
-E_p^hL^hu^h + cu^h + \gamma(I^h-E_q^h)u^h = f^h.
\end{equation}
We recall that the points in $B_{edge}$ do not belong to any of the interpolation stencils,
which means that the columns in $E_p^h$ corresponding to points in $B_{edge}$ are all zero columns,
so the rows in $L^h$ corresponding to points in $B_{edge}$ can vary arbitrarily while
$E_p^hL^h$ stays unchanged. Thus we see that no boundary condition for the discrete Laplace matrix
is needed, since the discrete Laplacian at any points in $B_{edge}$ does not contribute to
the closest point evaluation of $L^hu^h$.

As mentioned before, $\gamma$ is a parameter balancing the side and consistency conditions.
As pointed out in \cite{plain:icpm},
a typical choice of $\gamma$ is $\frac{2d}{\Delta x^2}$, where $d$ is the
dimension of the embedding space and $\Delta x$ is the grid size. In \cite{plain:icpm}, this choice of $\gamma$
is motivated by stabilizing the implicit time-stepping scheme for the in-surface heat equation.

The choice of $\gamma=\frac{2d}{\Delta x^2}$ is not obvious simply from (\ref{discrete}),
but it turns out to be a good choice there.
In Section~\ref{convergence-proof} we show that for the shifted Poisson's equation
on a closed curve embedded in $\mathbb{R}^2$, the coefficient matrix corresponding to (\ref{discrete})
is diagonally dominant by rows when we choose $\gamma=\frac{2d}{\Delta x^2}$, $p=1$, $q=3$,
and $L^h$ as the Laplacian matrix from the standard centered finite difference scheme.

\subsection{Truncation error}
\label{truncation-error}
In this section, we show that the difference scheme (\ref{discrete}) is consistent for the embedding equation (\ref{embedding})
provided that $\tilde{u}$ is sufficiently differentiable.

We denote the restriction operator which evaluates the continuous function at the grid points by $R^h$,
and denote the coefficient matrix in the difference scheme (\ref{discrete}) by $A^h$, i.e.,
\begin{equation}
\label{what-is-A}
A^h = cI^h-M^h, \quad\text{where}\quad M^h = E_p^hL^h - \gamma(I^h-E_q^h).
\end{equation}
Then we prove a theorem showing the form of the truncation error.
\begin{thm}
\label{truncation-thm}
Assume that $\tilde{u}\in C^{k}(\mathbb{R}^d)$ is a solution of (\ref{embedding}), where $d$ is the dimension
of the embedding space and $k\geq \max\{4, p+3, q+1\}$ ($p$ and $q$ are the
polynomial interpolation degrees for the discrete closest point extension of $\Delta\tilde{u}$ and $\tilde{u}$
respectively). Let $L^h$ be a second-order (in the $\infty$-norm) approximation to
the Cartesian Laplacian $\Delta$, and $\gamma=O(\frac{1}{\Delta x^2})$. We have the
following form of the truncation error,\footnote{If $\bm{v}$ is a column vector, then by $\bm{v}=O(\Delta x^k)$
we mean that $\|\bm{v}\|_{\infty}\leq C\Delta x^k$, where $k$ and $C$ are positive constants.}
\begin{equation}\label{truncation}
A^hR^h\tilde{u}(x)-R^h\tilde{f}(x)=O(\Delta x^{p+1}) + O(\Delta x^{q-1}) + O(\Delta x^2),
\end{equation}
where $A^h$ is as given in (\ref{what-is-A}), and $R^h$ is the restriction operator mentioned before.
\end{thm}

{\it Proof.} Let $w(x)=\Delta\tilde{u}(x)$, $x\in B(\mathcal{S})$; then $w(x)\in C^{k-2}(\mathbb{R}^d)$.
As mentioned before, since $\tilde{u}$ is constant along directions normal to $\mathcal{S}$,
$w(y)=\tilde{f}(y)$ for any $y\in\mathcal{S}$. We evaluate $w$ at the closest points on
$\mathcal{S}$ through degree $p$ interpolation, which will introduce an interpolation error.
Since $k\geq p+3$, $k-2\geq p+1$, we have the following standard error estimate for degree
$p$ polynomial interpolation,
\begin{equation}
\label{normal-rhs}
E_p^hR^hw(x) - R^h\tilde{f}(x) = O(\Delta x^{p+1}).
\end{equation}
Furthermore, since $\tilde{u}\in C^k(\mathbb{R}^d)$, $k\geq 4$,
the Laplacian matrix $L^h$ is a second-order approximation to $\Delta$,
we have that
\begin{equation}
\label{pde}
R^hw(x) = L^hR^h\tilde{u}(x) + O(\Delta x^2).
\end{equation}
Plugging (\ref{pde}) into (\ref{normal-rhs}), we get
\begin{displaymath}
E_p^h(L^hR^h\tilde{u}(x) + O(\Delta x^2) ) - R^h\tilde{f}(x) = O(\Delta x^{p+1}),
\end{displaymath}
Since we are using local polynomial interpolation, $\|E_p^h\|_{\infty}$ is bounded above by some constant
independent of $\Delta x$, so $E_p^hO(\Delta x^2)=O(\Delta x^2)$; and thus we get
\begin{equation}
\label{pde-final}
E_p^hL^hR^h\tilde{u}(x) - R^h\tilde{f}(x) = O(\Delta x^{p+1}) + O(\Delta x^2).
\end{equation}
Finally, we estimate the error of evaluating $\tilde{u}(\text{cp}(x))$ by degree $q$ interpolation.
Because $\tilde{u}\in C^{k}(\mathbb{R}^d)$ and $k\geq q+1$, once again we have the standard interpolation
error estimation,
\begin{equation}
\label{normal-u}
(I^h-E_q^h)R^h\tilde{u}(x) = O(\Delta x^{q+1}).
\end{equation}
Multiplying equation (\ref{normal-u}) by $\gamma=O(\frac{1}{\Delta x^2})$ and adding the results with equation (\ref{pde-final}),
we have that equation (\ref{truncation}) holds. $\hfill\blacksquare$ 

\begin{remark}
\label{remark-trunc}
According to (\ref{truncation}), if we pick $p=1$, $q=3$, we will get a difference scheme for the
embedding equation (\ref{embedding}) with second-order truncation error in the $\infty$-norm. 
Contrasting this discretization with the method of Macdonald--Brandman--Ruuth \cite{plain:eigen,plain:icpm}, we note that discretizing (\ref{embedding-icpm}) to second-order consistency requires
cubic interpolations in \emph{both} terms.
Because lower degree interpolations have smaller stencils, a discretization of (\ref{embedding}) will yield a coefficient matrix
that is sparser than that of (\ref{embedding-icpm}).
For example, for a curve in 2D, our new approach is roughly 34\% sparser and 50\% for a surface in 3D.
\end{remark}



\subsection{Convergence for a closed curve in $2$-D}
\label{convergence-proof}
In this section, we aim to prove convergence of the difference scheme (\ref{discrete}) in the restricted case of a closed curve in 2D.
\begin{thm} \label{curve-convergence}
Suppose that we wish to solve the embedding equation (\ref{embedding}) on a smooth closed curve in $\mathbb{R}^2$,
and the solution $\tilde{u}$ of the embedding equation
satisfies the assumptions in Theorem \ref{truncation-thm}. If in (\ref{discrete}) we choose
$\gamma=\frac{4}{\Delta x^2}$, $p=1$, $q=3$, and $L^h$ as the Laplacian matrix from
the standard $\frac{1}{\Delta x^2}\{-4, 1, 1, 1, 1\}$ finite difference stencil,
then the difference scheme (\ref{discrete}) is second-order convergent in the $\infty$-norm.
\end{thm}

If the embedding space is $\mathbb{R}^2$, when picking $\gamma=\frac{4}{\Delta x^2}$, and $p=1$, $q=3$,
the coefficient matrix becomes
\begin{equation}
\label{what-is-A-2d}
A^h = cI^h-M^h, \quad\text{where}\quad M^h = E_1^hL^h - \frac{4}{\Delta x^2}(I^h-E_3^h).
\end{equation}
We prove Theorem~\ref{curve-convergence} using the following theorem and two propositions.
\begin{thm}
\label{mmatrix}
If the embedding space is $\mathbb{R}^2$, then $A^h$ defined by (\ref{what-is-A-2d})
(i.e., $\gamma=\frac{4}{\Delta x^2}$, $p=1$, $q=3$, and $L^h$ is the Laplacian matrix from
the standard $\frac{1}{\Delta x^2}\{-4, 1, 1, 1, 1\}$ finite difference stencil) is an M-Matrix,
with positive diagonal entries $a_{kk}>0$,
non-positive off-diagonal entries $a_{kj}\leq 0$, $k\not=j$,
and each row sums up to $c$:
\begin{equation}
a_{kk} + \sum_{k\not=j}a_{kj} = c.
\end{equation}
\end{thm}

\begin{prop}(Varah 1975 \emph{\cite{plain:Varah}})
\label{inverse-bound}
Assume that $A$ is diagonally dominant by rows and set $c =\emph{min}_k(|a_{kk}|-\sum_{j\not=k}|a_{kj}|)$.
Then $\|A^{-1}\|_{\infty}\leq 1/c$.
\end{prop}

\begin{prop}
\label{convergence}
If $\|(A^h)^{-1}\|_{\infty}\leq 1/c$, then the difference scheme  (\ref{discrete})
is second-order convergent in the $\infty$-norm provided that the scheme has second-order truncation error
in the $\infty$-norm.
\end{prop}

Theorem \ref{mmatrix} indicates that the coefficient matrix $A^h$ satisfies the assumption of
Proposition \ref{inverse-bound}, which means that the claim of Proposition \ref{inverse-bound} holds:
$\|(A^h)^{-1}\|_{\infty}\leq 1/c$; and this is in turn the assumption of Proposition \ref{convergence},
so the claim of Proposition \ref{convergence} holds. Putting the theorem and propositions together,
Theorem \ref{curve-convergence} holds.


\paragraph{Proof of Proposition~\ref{convergence}}
Choosing $p=1$, $q=3$, the truncation error estimate (\ref{truncation}) becomes
$A^hR^h\tilde{u}(x) - R^h\tilde{f}(x) = O(\Delta x^2)$.
From the discrete linear system (\ref{discrete}) we have that
$A^hu^h - R^h\tilde{f}(x) = 0$.
Subtracting these two results gives
$A^h(u^h-R^h\tilde{u}(x)) = O(\Delta x^2)$.
Because the $\infty$-norm of $(A^h)^{-1}$ is bounded by a positive constant $1/c$, which
is independent of $\Delta x$, we get
$u^h-R^h\tilde{u}(x) = (A^h)^{-1}O(\Delta x^2) = O(\Delta x^2)$,
i.e., the difference scheme (\ref{discrete}) is second-order convergent in the $\infty$-norm. $\hfill\blacksquare$
\bigskip

\paragraph{Proof of Theorem~\ref{mmatrix}}
The proof contains two steps.

\paragraph{Step 1} We first show that {\it each row of $M^h$ sums up to zero}.

Taking all $u_i$ to be $1$ in equation (\ref{pointwise}), the left-hand side minus $c$ becomes
\begin{equation}
\label{row-sum}
-\sum_{x_j\in\text{I}_p(\text{cp}(x_i))} \omega_j
\Bigg(\sum_{x_k\in\text{Diff}(x_j)} l_k\Bigg)
+ \gamma\Bigg(1- \sum_{x_m\in\text{I}_q(\text{cp}(x_i))}\bar{\omega}_m\Bigg).
\end{equation}
Consistency of the finite difference scheme implies that
$\sum_{x_k\in\text{Diff}(x_j)} l_k = 0, \quad \text{for any} \quad x_j\in L$.
By the construction of the interpolation weights,
$\sum_{x_m\in\text{I}_q(\text{cp}(x_i))}\bar{\omega}_m = 1$.
So (\ref{row-sum}) equals $0$, i.e., each row of $M^h$ sums up to zero.

\paragraph{Step 2} We then prove that {\it all the off-diagonal entries of $M^h$ are non-negative},
thus all the diagonal entries are non-positive since each row of $M^h$ sums up to zero.
For convenience (not necessity) we first rearrange $M^h$ as
\begin{displaymath}
M^h = E_1^h(L^h+\frac{4}{\Delta x^2}I^h) + \frac{4}{\Delta x^2}(E_3^h-E_1^h) - \frac{4}{\Delta x^2}I^h.
\end{displaymath}

\begin{figure}[htbp]
\begin{center}
\subfigure[Relative position of the closest point $\diamond$ in the stencil]{
\label{stencil:a}
\begin{tikzpicture}[scale=3]
\draw (0.83,0) arc (0:90:0.83cm);
\draw (0.83,0) -- (0.83,-0.1);
\node [right] at (0.83,-0.1) {$\mathcal{S}$};
\draw (0,0.83) -- (-0.1,0.83);
\foreach \i in {0,...,3}
{
    \foreach \j in {0,...,3}
	{
	\draw plot[mark = *, mark size = 0.4 pt] coordinates{(1/3*\i,1/3*\j)};	
	}	
}
\draw plot[mark = diamond, mark size = 0.8pt] coordinates{(0.586898628,0.586898628)};
\foreach \i in {0,...,3}
{
    \draw (1/3*\i,0) -- (1/3*\i,1);	
}
\foreach \i in {0,...,3}
{
    \draw (0, 1/3*\i) -- (1, 1/3*\i);	
}
\node [left] at (0,0.5) {$\Delta x$};
\node [below] at (0.5,0) {$\Delta x$};
\draw [style=dashed] (0.586898628,0.586898628) -- (0.3333333333,0.586898628)
  node [below] at (0.460115981,0.7) {$a$};
\draw [style=dashed] (0.586898628,0.586898628) -- (0.586898628,0.3333333333)
  node [left] at (0.7,0.460115981) {$b$};
\end{tikzpicture}
}
\hspace{0.5cm}
\subfigure[Local index for points in the stencil]{
\label{stencil:b}
\begin{tikzpicture}[scale=3]
\tikzstyle{every node}=[font=\small]
\draw (0.83,0) arc (0:90:0.83cm);
\draw (0.83,0) -- (0.83,-0.1);
\draw (0,0.83) -- (-0.1,0.83);
\foreach \i in {0,...,3}
{
    \foreach \j in {0,...,3}
	{
	\draw plot[mark = *, mark size = 0.4 pt] coordinates{(1/3*\i,1/3*\j)};	
	}	
}
\draw plot[mark = diamond, mark size = 0.8pt] coordinates{(0.586898628,0.586898628)};
\foreach \i in {0,...,3}
{
    \draw (1/3*\i,0) -- (1/3*\i,1);	
}
\foreach \i in {0,...,3}
{
    \draw (0, 1/3*\i) -- (1, 1/3*\i);	
}
\foreach \i in {0,...,3}
{
    \foreach \j in {0,...,3}
	{
	\node [below] at (1/3*\i,1/3*\j) {\j\i};	
	}	
}
\end{tikzpicture}
}
\hspace{0.5cm}
\subfigure[Signs of the weights of bi-cubic interpolation for
the closest point $\diamond$]{
\label{stencil:c}
\begin{tikzpicture}[scale=3]
\tikzstyle{every node}=[font=\small]
\foreach \i in {0,...,3}
{
    \foreach \j in {0,...,3}
	{
	\draw plot[mark = *, mark size = 0.4 pt] coordinates{(1/3*\i,1/3*\j)};	
	}	
}
\draw plot[mark = diamond, mark size = 0.8pt] coordinates{(0.586898628,0.586898628)};
\foreach \i in {0,...,3}
{
    \draw (1/3*\i,0) -- (1/3*\i,1);	
}
\foreach \i in {0,...,3}
{
    \draw (0, 1/3*\i) -- (1, 1/3*\i);	
}
\node [below right] at (0,0) {$+$};
\node [below right] at (0.3333333333,0) {$-$};
\node [below right] at (0.6666666667,0) {$-$};
\node [below right] at (1,0) {$+$};
\node [below right] at (0,0.3333333333) {$-$};
\node [below right] at (0.3333333333,0.3333333333) {$+$};
\node [below right] at (0.6666666667,0.3333333333) {$+$};
\node [below right] at (1,0.3333333333) {$-$};
\node [below right] at (0,0.6666666667) {$-$};
\node [below right] at (0.3333333333,0.6666666667) {$+$};
\node [below right] at (0.6666666667,0.6666666667) {$+$};
\node [below right] at (1,0.66666666667) {$-$};
\node [below right] at (0,1) {$+$};
\node [below right] at (0.3333333333,1) {$-$};
\node [below right] at (0.6666666667,1) {$-$};
\node [below right] at (1,1) {$+$};
\end{tikzpicture}
}
\end{center}
\caption[Local information for a cubic interpolation stencil]
{Information for a $4\times 4$ stencil surrounding the closest point {$\diamond$}. }
\label{stencil}
\end{figure}
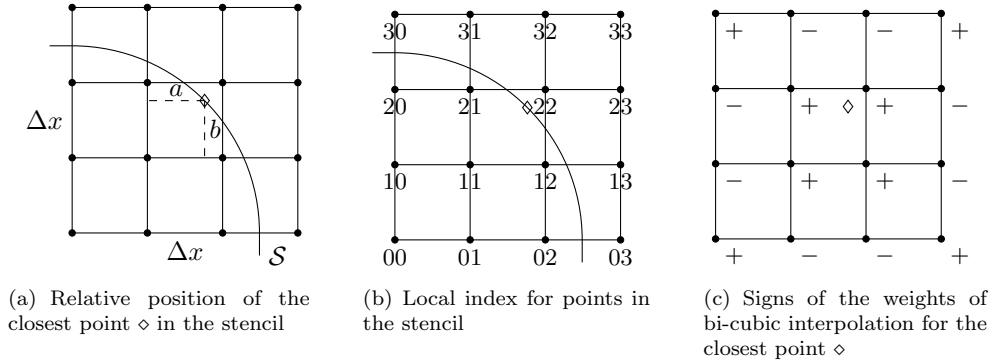

We then investigate the result of multiplying the $i$-th row of $M^h$ by the column vector $u^h$.
Figure~\ref{stencil}(a) shows the position of
$\text{cp}(\bm{x}_i)$ (indicated by $\diamond$). $\text{cp}(\bm{x}_i)$ is in the middle
block of the bi-cubic interpolation stencil, and its relative position in the stencil is as shown
in Figure~\ref{stencil}(a).
For clarity of notation, we number the points in the bi-cubic interpolation stencil by a local index
as shown in Figure~\ref{stencil}(b). Denote the points in the stencil with respect to the local index
by $X_{lm}$ ($0\leq l,m\leq 3$),
the corresponding values at these points by $U_{lm}$ ($0\leq l,m\leq 3$), the bi-cubic
interpolation weights at points $X_{lm}$ by $\bar{\omega}_{lm}$ ($0\leq l,m\leq 3$), and the bi-linear
interpolation weights at points $X_{nk}$ by $\omega_{nk}$ ($0\leq n,k\leq 1$).
The $i$-th entry of $M^hu^h$ becomes
\begin{gather}
\omega_{11}(U_{10}+U_{01}+U_{12}+U_{21})/\Delta x^2 \nonumber\\
+\omega_{12}(U_{11}+U_{02}+U_{13}+U_{22})/\Delta x^2 \nonumber\\
+\omega_{21}(U_{20}+U_{11}+U_{22}+U_{31})/\Delta x^2 \nonumber\\
+\omega_{22}(U_{21}+U_{12}+U_{23}+U_{32})/\Delta x^2 \nonumber\\
+4\Big(\sum_{0\leq l,m\leq 3}\bar{\omega}_{lm}U_{lm}-\sum_{1\leq l,m \leq 2}\omega_{lm}U_{lm}\Big)/\Delta x^2 - 4u^h_i/\Delta x^2.
\end{gather}
 By the notation mentioned before, $U_{lm}$ ($0\leq l,m\leq 3$) are just the values of $u^h$
 at the points in the bi-cubic interpolation stencil written in the local index.
 The point $\bm{x}_i$ might not be contained in the interpolation stencil, so we still write $u^h_i$
 in its global index.
As $l,m$ go from $1$ to $3$, $X_{lm}$ might or might not be $\bm{x}_i$, depending on whether the grid
point $\bm{x}_i$ is in the interpolation stencil of $\text{cp}(\bm{x}_i)$ or not. In either case the
negative coefficient of $u^h_i$ only contributes to the diagonal entry in the corresponding row of the matrix $M^h$,
but the coefficient does not contribute to the off-diagonal entries of $M^h$.
So we only need to show that {\it all the coefficients of $U_{lm}$, excluding the potential
contribution from the negative coefficient of $u^h_i$, are non-negative for all $0\leq l,m\leq 3$.
}

If $\text{cp}(\bm{x}_i)$ coincides with some grid point $\bm{x}_j$, then all the interpolation weights
vanish except that the weight at $\bm{x}_j$ is $1$, so it is obvious that the above claim is true. Now we
only consider cases where $\text{cp}(\bm{x}_i)$ does not coincide with any grid point, which is the
case illustrated in Figure~\ref{stencil}.

In the rest of the proof, it might be helpful to keep in mind the signs of the bi-cubic interpolation
weights $\bar{\omega}_{lm}$ ($0\leq l,m \leq 3$) for $\text{cp}(\bm{x}_i)$  at each point of the stencil.
Figure~\ref{stencil}(c) shows the signs of the bi-cubic
interpolation weights for $\text{cp}(\bm{x}_i)$ (indicated by $\diamond$) which is in the
middle block of the bi-cubic interpolation stencil. These signs can be obtained by straightforward computations,
and they do not depend on the precise location of $\text{cp}(\bm{x}_i)$.
Since we interpolate in a dimension-by-dimension fashion,
the weights at each of the points are obtained by multiplication of the corresponding weights
in the $x$ and $y$ directions. For instance, if we want to compute the bi-cubic interpolation
weight at the lower-left corner point $\bar{\omega}_{00}$,
we first compute the weight in the $x$ direction,
$\bar{\omega}_{00}^x = -\frac{a(\Delta x-a)(2\Delta x-a)}{6\Delta x^3}$,
and the weight in the $y$ direction,
$\bar{\omega}_{00}^y = -\frac{b(\Delta x-b)(2\Delta x-b)}{6\Delta x^3}$,
and then we multiply the weights in the two directions to get
$\bar{\omega}_{00} = \bar{\omega}_{00}^x \bar{\omega}_{00}^y$.
 Here we have used the superscripts $x$ and $y$ in the notations $\bar{\omega}_{00}^x$ and $\bar{\omega}_{00}^y$
 to indicate the weights in the $x$ and $y$ directions, respectively. We shall
 use this notation in the following proof.
The bi-cubic weights at other points, and the bilinear weights, can be computed in the same way. Some of
their values will be presented in the rest of the proof when necessary.

Now let us finish the proof by calculating the coefficients of $U_{lm}$.
There are three cases for the coefficients of $U_{lm}$ as $l,m$ go from $1$ to $3$. We only need to
verify that for each of the three cases, the coefficients of $U_{lm}$ are positive.
\begin{itemize}
\item[1.] Four corner points:
 $(l,m) = (0,0), (0,3), (3,0), (3,3)$.\\
 The only contribution for the row entries corresponding to
 these four points are the bi-cubic interpolation weights
 $\bar{\omega}_{00}$, $\bar{\omega}_{03}$, $\bar{\omega}_{30}$, $\bar{\omega}_{33}$,
 which happen to be positive in $2$-D.
\item[2.] Four center points:
 $(l,m) = (1,1), (1,2), (2,1), (2,2)$.\\
 There are some positive contributions from the linear combinations of the weights in
 the centered difference scheme;
 so we only need to show at each of these four center points,
 the positive weight of the bi-cubic interpolation
 minus the bi-linear interpolation weight is larger than zero.
 For the sake of symmetry, we only show that
 $\bar{\omega}_{11}-\omega_{11}>0$,
 $\text{i.e.}, \bar{\omega}_{11}^x\bar{\omega}_{11}^y-\omega_{11}^x\omega_{11}^y>0$.
 Since $\bar{\omega}_{11}^x$, $\bar{\omega}_{11}^y$, $\omega_{11}^x$, and $\omega_{11}^y$ are positive,
 it is enough to show that
 \begin{equation}
 \label{less_tedious}
 \bar{\omega}_{11}^x>\omega_{11}^x \quad \text{and} \quad \bar{\omega}_{11}^y>\omega_{11}^y.
 \end{equation}
 By straightforward computation,
 $\bar{\omega}_{11}^x=\frac{(\Delta x+a)(\Delta x-a)(2\Delta x-a)}{2\Delta x^3}$, 
 $\omega_{11}^x = \frac{\Delta x-a}{\Delta x}$, 
 $\bar{\omega}_{11}^y=\frac{(\Delta x+b)(\Delta x-b)(2\Delta x-b)}{2\Delta x^3}$, 
 $\omega_{11}^y = \frac{\Delta x-b}{\Delta x}$. 
 It is then straightforward to verify that (\ref{less_tedious}) holds for any $0<a<\Delta x$ and $0<b<\Delta x$.
\item[3.] Eight edge points:
 $(l,m) = (1,0),(2,0),(0,1),(0,2),(1,3),(2,3),(3,1),(3,2)$.\\
 Again by symmetry we only need to show that the row entry corresponding to $(l,m)=(1,0)$ is positive, i.e.,
$\frac{\omega_{11}}{\Delta x^2} + \frac{4\bar{\omega}_{10}}{\Delta x^2} > 0$,
 where $\omega_{11}>0$, but $\bar{\omega}_{10}<0$.
 So we need to show that $|\omega_{11}|>4|\bar{\omega}_{10}|$, i.e.,
$|\omega_{11}^x\omega_{11}^y| > 4|\bar{\omega}_{10}^x\bar{\omega}_{10}^y|$.
 In fact we can show that
 \begin{equation}
 \label{tedious}
 \omega_{11}^x>6|\bar{\omega}_{10}^x|, \quad \text{and} \quad \omega_{11}^y\geq\frac{8}{9}|\bar{\omega}_{10}^y|.
 \end{equation}
 Again by straightforward computation,
 $\omega_{11}^x = \frac{\Delta x-a}{\Delta x}$, 
 $|\bar{\omega}_{10}^x| = \frac{a(\Delta x-a)(2\Delta x-a)}{6\Delta x^3}$, 
 $\omega_{11}^y = \frac{\Delta x-b}{\Delta x}$, 
 $|\bar{\omega}_{10}^y|=\frac{(\Delta x+b)(\Delta x-b)(2\Delta x-b)}{2\Delta x^3}$. 
 Again one can straightforwardly verify that (\ref{tedious}) holds for any $0<a<\Delta x$ and $0<b<\Delta x$.
\end{itemize}

Thus, Theorem \ref{mmatrix} holds:
the diagonal entries of $A^h=cI^h-M^h$ are all positive, the off-diagonal entries are all non-positive,
and each row sums up to $c$. $\hfill\blacksquare$

\begin{remark}
The M-Matrix result of Theorem \ref{mmatrix}, at least for our particular parameter choices, does not hold in 3-D.
The barycentric interpolation weights at the corners of the interpolation stencils are negative. A different convergence
proof would be required.
However, in practice we still observe second-order convergence in 3-D (Section~\ref{numerical-tests-multigrid}).
\end{remark}


\subsection{Numerical examples in $\mathbb{R}^2$}\label{numerics-2d}


\subsubsection{Shifted Poisson's equation on the unit circle}\label{poisson-circle}
We solve the shifted Poisson's equation $-\Delta_{\mathcal{S}}u + u = f$ on the unit circle $\mathcal{S}$;
the exact solution is chosen as $\sin(\theta)+\sin(12\theta)$, and we compute the right-hand side function $f$
from the exact solution. The discrete right-hand side is then obtained by evaluating the closest point extension of $f$
at the grid points in the computational band.
We construct $M^h=E^h_1L^h-\frac{4}{\Delta x^2}(I^h-E^h_3)$, and the coefficient matrix
is $A^h=I^h-M^h$.

Since the embedding equation is a $2$D problem,
sparse direct solvers (e.g.\ Matlab's backslash) work well in solving the resulting linear system.
After solving the linear system, we obtain the discrete solution defined on the banded grid points.
To calculate the errors, we then place many sample points on $\mathcal{S}$, estimate the numerical solution at these points
by cubic interpolation of the values at the surrounding grid points, and compute the $\infty$-norm of the errors
at the sample points on $\mathcal{S}$. Table~\ref{convergence-circle-table} shows second order convergence
of the relative errors in the $\infty$-norm. Alternatively, the exact closest point extension of the exact solution can be performed
onto the computational grid, and the errors can be directly calculated at the grid points; similar results are observed in this case.

\begin{table}[htbp]
\caption{\label{convergence-circle-table}
Convergence study for the shifted Poisson's equation on a unit circle showing 2nd-order convergence.}
\centering
\begin{small}
\begin{tabular}{cccccccc}
\hline
$\Delta x$                  & 0.1     & 0.05    &   0.025   & 0.0125  & 0.00625 & 0.003125 & 0.0015625\\
\hline
\rule[-1.6ex]{0pt}{4.6ex}  
$\frac{\|u-u^h\|_{\infty}}{\|u\|_{\infty}}$ & 6.51e-2 & 1.85e-2 &   4.76e-3 & 1.19e-3 & 2.97e-4 & 7.38e-5 & 1.85e-5\\
\rule[-1.7ex]{0pt}{4.6ex}  
$\log_2\frac{\|u-u^{2h}\|_{\infty}}{\|u-u^{h}\|_{\infty}}$ & & 1.81 & 1.95 & 2.00 & 2.00 & 2.01 & 2.00\\
\hline
\end{tabular}
\end{small}
\end{table}

\subsubsection{Shifted Poisson's equation on a bean-shaped curve} \label{bean-backslash}
\begin{figure}[htbp]
\begin{center}
\includegraphics[width=0.47\textwidth]{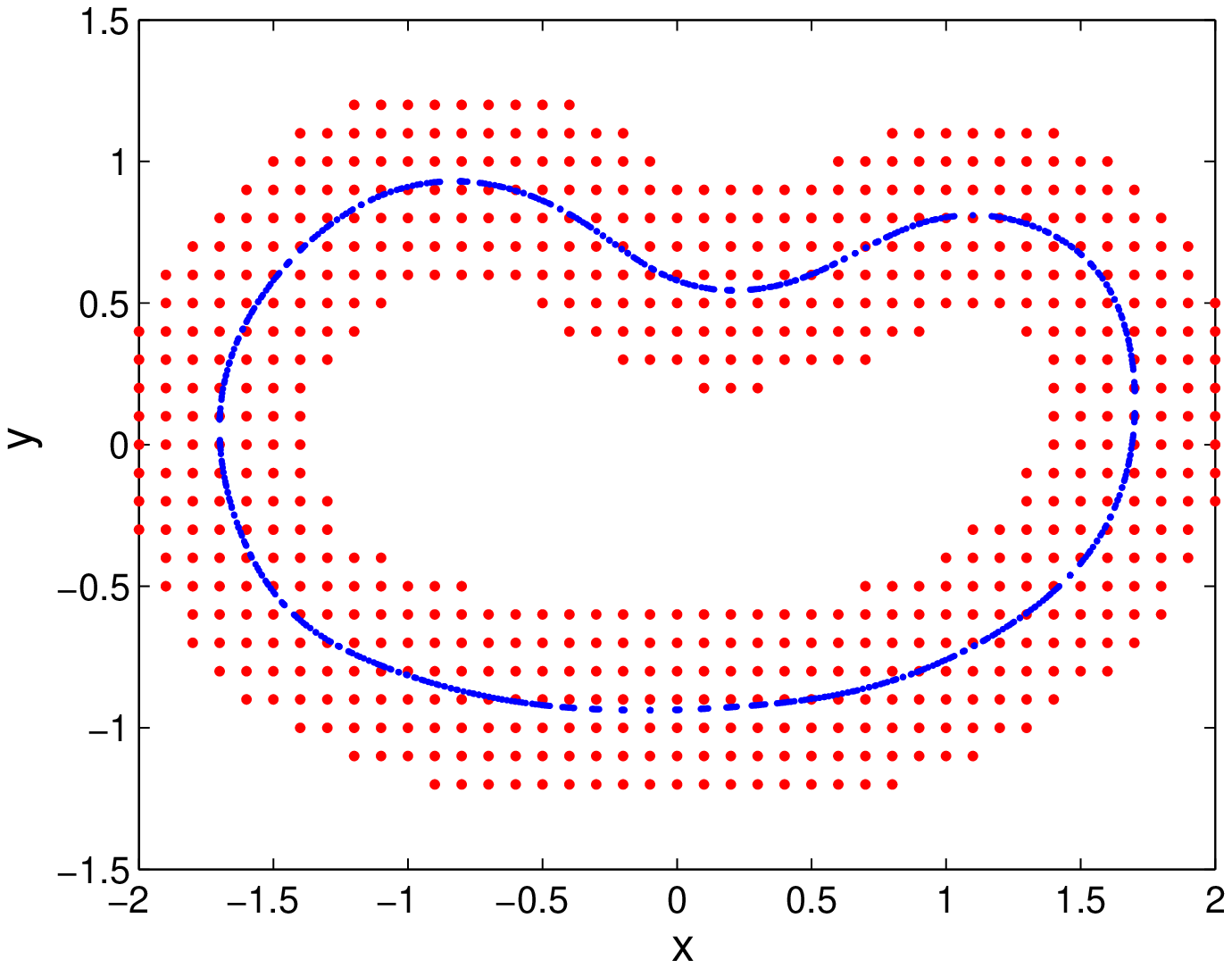}
\includegraphics[width=0.47\textwidth]{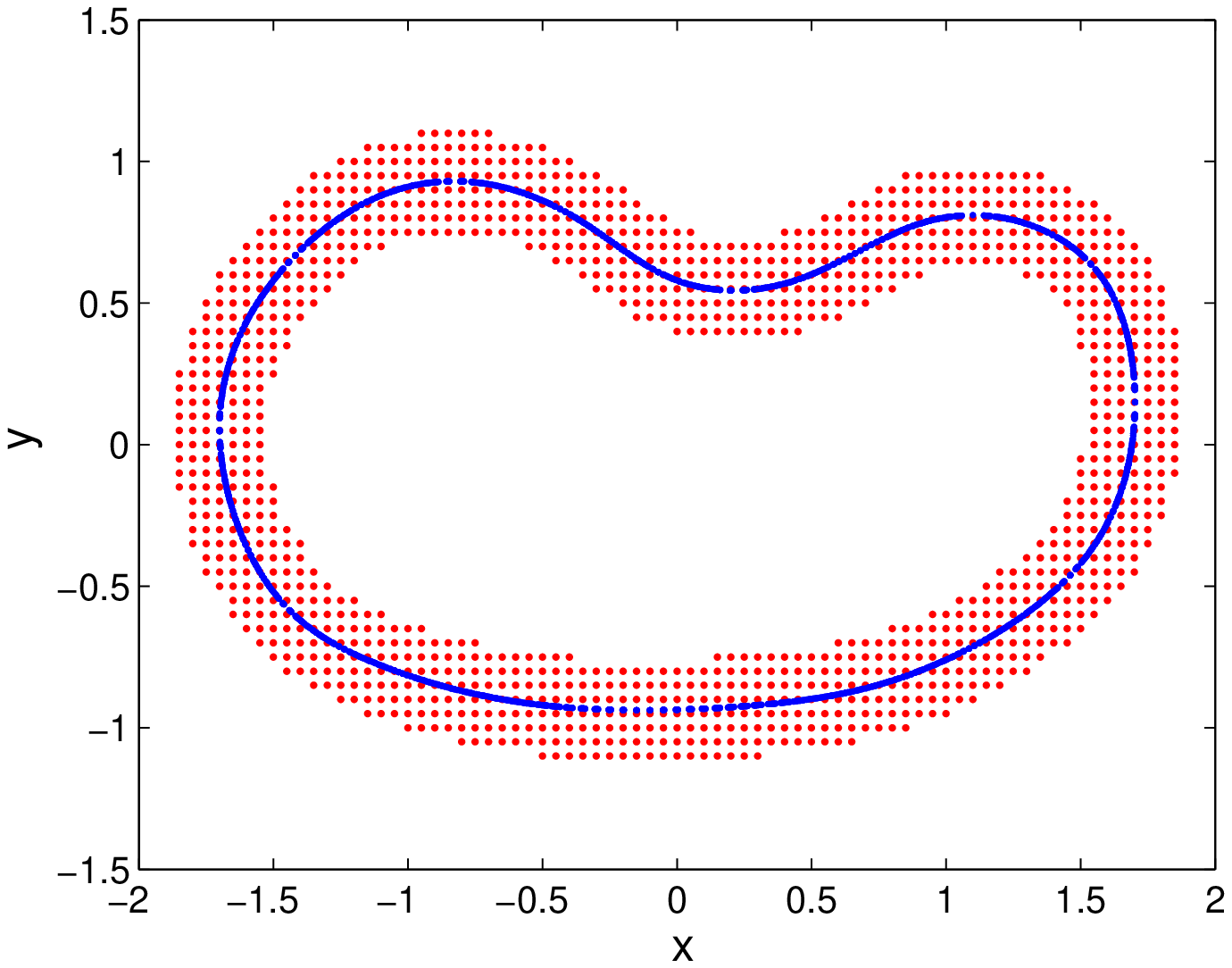}
\caption[Computational grids for a bean-shaped curve]{\label{bean-grid}
Computational grids for the Closest Point Method
on a bean-shaped curve with $\Delta x=0.1$ (left) and $\Delta x=0.05$ (right).}
\end{center}
\end{figure}

Again we solve the shifted Poisson's equation $-\Delta_{\mathcal{S}}u + u = f$, but this time on a bean-shaped curve (see Figure~\ref{bean-grid}).

This bean-shaped curve is constructed from a sequence of control points in $\mathbb{R}^2$;
the Matlab function `cscvn' \cite{plain:matlab-cscvn} is used
to get a parametric cubic spline curve passing through these points.
This results in $(x(t),y(t))$ parameterized by $t\in[0,T]$.
We obtain a closest point representation from this
parametrization as follows: we fix a grid point, then construct the
squared distance as a function of the parameter $t$.
Next we minimize that function (over the parameter $t$) using
Newton's method \cite{plain:boyd}.
This is repeated for each grid point in the computational band.

We construct the exact solution in terms of the parametrization, $u(t)=\sin(\frac{2\pi}{T}t)+\sin(\frac{20\pi}{T}t)$.
The Laplace--Beltrami operator of $u$ in terms of this parametrization is
$\Delta_\mathcal{S} u(t) = \frac{1}{\sqrt{\dot{x}(t)+\dot{y}(t)}}
                          \frac{\text{d}}{\text{d}t}
						  \big(  \frac{1}{\sqrt{\dot{x}(t)+\dot{y}(t)}} \frac{\text{d}u(t)}{\text{d}t} \big)$.
The right-hand side function is $f(t)=-\Delta_\mathcal{S}u(t)+u(t)$.
Again we construct the discrete right-hand side and coefficient matrix, solve the resulting linear system to get the numerical
solution, and compute the $\infty$-norm of the error on the bean-shaped curve. 
Table~\ref{convergence-bean-table} shows second order convergence results.
\begin{table}[htbp]
\caption[Shifted Poisson's equation on a bean-shaped curve]{\label{convergence-bean-table}
Convergence study for the shifted Poisson's equation on a bean-shaped curve showing 2nd-order convergence.}
\centering
\begin{small}
\begin{tabular}{cccccccc}
\hline
$\Delta x$                  & 0.1     & 0.05    &   0.025   & 0.0125  & 0.00625 & 0.003125 & 0.0015625\\
\hline
\rule[-1.6ex]{0pt}{4.6ex}  
$\frac{\|u-u^h\|_{\infty}}{\|u\|_{\infty}}$ & 6.19e-2 & 1.17e-2 &  3.59e-3 & 8.77e-4 & 2.35e-4 & 5.32e-5 & 1.38e-5\\
\rule[-1.7ex]{0pt}{4.6ex}  
$\log_2\frac{\|u-u^{2h}\|_{\infty}}{\|u-u^{h}\|_{\infty}}$ & & 2.40 & 1.71 & 2.03 & 1.90 & 2.14 & 1.95\\
\hline
\end{tabular}
\end{small}
\end{table}

\section{V-Cycle Multigrid solver}
\label{section-v-cycle}
In Section~\ref{numerics-2d} we studied some numerical examples on curves in $2$D, where the resulting
linear systems were solved by sparse direct solvers (e.g., Matlab's backslash).
However, for problems where the embedding space is $\mathbb{R}^3$ or higher,
iterative solvers are usually more desirable.
Since the geometric multigrid method (see e.g., \cite{plain:briggs,plain:hackbusch})
is quite effective in solving linear systems arising from discretizations of elliptic operators,
we combine it with the Closest Point Method to solve the linear system arising from the embedding equation
corresponding to the surface elliptic PDEs.

\subsection{The Multigrid ``smoother'': Jacobi iteration}
One of the key ingredients of a Multigrid algorithm is the ``smoother'' \cite{plain:briggs,plain:hackbusch},
which damps out high-frequency errors and thus smooths the solution at each grid level.
One good choice of ``smoother'' is the Jacobi iteration: suppose we want to solve $A^hu^h=f^h$,
starting from some initial $u^h$, we update the iterates by
\begin{equation}
\label{standard-relax}
u^h := u^h + \text{diag}(A^h)^{-1}(f^h - A^hu^h).
\end{equation}
One possible approach is to simply use this standard Jacobi iteration
on the discretization (\ref{what-is-A}) of our embedding equation (\ref{embedding}).
We shall discuss this standard relaxation strategy further in Section~\ref{section-standard-relax}.


An alternative strategy, and the one we propose here is based instead
on considering the earlier constrained system (\ref{pde-constraint-system}):
we want to ensure the solution satisfies the constraint
(\ref{constraint}) so we re-extend the solution following each Jacobi
iteration.


\subsection{Ruuth--Merriman Jacobi iteration}
\label{section-r-m}
The idea of this new smoothing strategy stems from the original Ruuth--Merriman approach for time-dependent surface PDEs \cite{plain:RuuthMerriman}.
In the continuous setting, if the solution $\tilde{u}$ is constant along the normal direction to the surface,
then the Cartesian Laplacian of $\tilde{u}$ agrees with the Laplace--Beltrami function of $\tilde{u}$ on the surface.
So roughly speaking $L^hu^h$ is consistent with $\Delta_{\mathcal{S}}\tilde{u}$ \emph{on the surface}, provided that
$u^h$ is constant along the normals to the surface.
Under this side condition, we can replace $M^h$ with $L^h$ in one iteration.
Let $\tilde{A}^h=cI^h-L^h$ be the shifted Cartesian Laplacian matrix.
Starting from an initial $u^h$ which is constant along the normal direction to the surface,
in one iteration we do the standard Jacobi iteration with respect to the $\tilde{A}^h=cI^h-L^h$ (instead of $A^h=cI^h-M^h$),
and then immediately do a closest point extension to ensure the discrete solution
is still constant along the normal direction to the surface:
\begin{subequations}
\begin{align}
\label{alg-relax}
u^h &:= u^h + \text{diag}(\tilde{A}^h)^{-1}(f^h - \tilde{A}^hu^h), \\
\label{alg-extend}
u^h &:= E^hu^h.
\end{align}
\end{subequations}

We remark that interpolation $E^h$ with sufficiently high degree should be used in order to avoid
introducing interpolation errors which dominate over the errors from the discrete Laplacian $L^h$.
In particular, we need to use cubic interpolation to build $E^h$
if we use a second-order difference scheme to build $L^h$.


\subsection{Restriction and prolongation operators}
\label{section-r-p}
In a Multigrid algorithm, one also needs to build \emph{restriction} operators which restrict data from fine grid
to coarse grid and \emph{prolongation} operators which prolong data from coarse grid to fine grid.

For practical efficiency, we would like to perform the computation in narrow bands with bandwidths shrinking
proportional to the grid sizes.
This makes the coarse grid band wider than the fine grid band (see Figure~\ref{rp}).
Some coarse grid points are out of the range of the fine grid band, which means there are no fine grid points
surrounding those ``outer'' coarse grid points; in those cases the standard restriction strategies such as ``injection'' and
``full weighting'' \cite{plain:briggs} do not work. We overcome this problem with the help of the closest point extension.
When restricting the values of points on the fine grid to the coarse grid, we assign the value
of each coarse grid point to be the value of its corresponding closest point, which
in turn is obtained by interpolating the values of fine grid points surrounding the
closest point. Naturally, we can construct the prolongation operators in the same
way:\footnote{An alternative approach is simply to use one of the standard nest-grid multigrid approaches (although we have not tested this here).}
take the value of each fine grid point to be the value of its corresponding closest point,
which is obtained through interpolation of values of coarse grid points.
Similarly to the closest point extension on a single grid, the restriction and prolongation operators
have corresponding matrix forms; we denote the restriction matrix by $E_{h}^{2h}$, and
the prolongation matrix by $E_{2h}^h$.
Figure~\ref{rp} illustrates the process of doing restriction and prolongation.

\begin{figure}[htbp]

\begin{center}
\begin{minipage}{6cm}
\begin{tikzpicture}[scale=3]
\draw[color=blue] plot[only marks, mark=o,mark size=0.6pt] file {coarse.dat}
  node [right] {coarse};

\draw[color=red] plot[only marks, mark=*,mark size=0.2pt] file {fine.dat}
  node [below] at (1,0)  {fine};

\draw (1,0) arc (0:90:1cm)
node [right]{};

\fill[red, opacity = 0.3] (0.9,0.3) rectangle(1.0,0.4)
  node [right]{};

\draw[style = dashed] (0.948683298050514,0.316227766016838) -- (1.2,0.4)
  node [right] {restriction};

\draw plot[mark = triangle, mark size = 0.4pt] coordinates{(0.948683298050514,0.316227766016838)}
  node [right] {};

\draw[style = dashed] (0.341743063086704,0.939793423488437) -- (0.4,1.1)
  node [above] at(0.2,1.2) {prolongation};

\draw plot[mark = triangle, mark size = 0.4pt] coordinates{(0.341743063086704,0.939793423488437)}
  node [right] {};

\fill[blue, opacity = 0.3] (0.2,0.8) rectangle(0.4,1.0)
  node [left] {};
\end{tikzpicture}
\end{minipage}
\begin{minipage}{6cm}
\caption[Illustration of Restriction and Prolongation]{\label{rp}
Illustration of Restriction and Prolongation. Blue circles indicate
the coarse grid points, red dots indicate the fine grid points.\\
Restriction: restrict quantities from fine grid to coarse grid using $E_h^{2h}$.\\
Prolongation: prolong quantities from coarse grid to fine grid using $E_{2h}^h$.}
\end{minipage}
\end{center}
\end{figure}

\subsection{Closest Point Method V-Cycle algorithm}
Armed with the smoothing strategy in Section~\ref{section-r-m} and the construction of
the restriction and prolongation operators in Section~\ref{section-r-p}, we
propose the Closest Point Method Multigrid algorithm.
In this paper, we use the simplest V-Cycle Multigrid scheme,
although the results could be extended to more general Multigrid schemes such as W-Cycle or the full Multigrid schemes.
We remark that this algorithm does not explicitly involve the side condition as in (\ref{embedding}),
but instead follows each iteration with a closest point extension as in the original Ruuth--Merriman approach \cite{plain:RuuthMerriman}.

\begin{alg}[Closest Point Method V-Cycle Algorithm]
\label{alg-vc-rm}
Denote by $B^h$ and $B^{2h}$ the fine and next coarser grids.
Denote by $E^h$ and $E^{2h}$ the extension matrices on the fine and next coarser grids (see Sec.~\ref{sec:diffscheme}).
Denote the shifted Cartesian Laplacian matrix, the solution and right-hand side on the fine grid by
$\tilde{A}^h=cI^h-L^h$, $u^h$, and ${f}^h$, and the corresponding variables on
the next coarser grid by $\tilde{A}^{2h}=cI^{2h}-L^{2h}$, $u^{2h}$, and ${f}^{2h}$.
Let ${E}_h^{2h}$ and ${E}_{2h}^h$ be the restriction and prolongation matrices (which are also extensions as in Sec.~\ref{section-r-p}).
We compute $u^h=V^h(u^h,f^h)$ by recursively calling:
\begin{itemize}
\item[1.]
Starting from some initial guess $u^h$ which is constant along the normals to the surface (e.g.\ $u^h=0^h$) on the finest grid,
do $\nu_1$ steps of Jacobi iteration with respect to $\tilde{A}^h$,
each step followed by a closest point extension:
\begin{eqnarray}
for\ i = 1:\nu_1,\ do & & \nonumber \\
u^h & := & u^h + \text{diag}(\tilde{A}^h)^{-1}({f}^h - \tilde{A}^hu^h); \nonumber \\
u^h & := & {E}^hu^h. \nonumber
\end{eqnarray}

\item[2.] Compute the residual and restrict it to the coarser grid,
          \begin{eqnarray}
          \label{residual-L}
		  r^h & := & {f}^h - \tilde{A}^h u^h, \\
          {f}^{2h} & := & {E}_h^{2h}r^h, \nonumber
          \end{eqnarray}
		  If $B^{2h} = \text{coarsest grid}$, solve $A^{2h} u^{2h} = {f}^{2h}$ directly; else
		  \begin{displaymath}
		  \begin{array}{rcl}
          u^{2h} & := & 0^{2h}, \\
          u^{2h} & := & V^{2h}(u^{2h},{f}^{2h}).
          \end{array}
          \end{displaymath}

\item[3.] Correct the fine grid solution by coarse grid solution:
\begin{displaymath}
u^h := u^h + {E}_{2h}^h u^{2h}.
\end{displaymath}

\item[4.] Do $\nu_2$ steps of Jacobi iteration with respect to $\tilde{A}^h$, each step followed by a closest point extension:
\begin{eqnarray}
for\ i = 1:\nu_2,\ do & & \nonumber \\
u^h & := & u^h + \text{diag}(\tilde{A}^h)^{-1}({f}^h - \tilde{A}^hu^h); \nonumber \\
u^h & := & {E}^hu^h. \nonumber
\end{eqnarray}
\end{itemize}
\end{alg}

\subsection{An alternative relaxation strategy}
\label{section-standard-relax}
In Step 1 and Step 4, instead of using the Ruuth--Merriman smoothing strategy (\ref{alg-relax}) and (\ref{alg-extend}),
one can simply use the iteration scheme (\ref{standard-relax}) which only depends on the linear system $A^hu^h=f^h$.
One can also easily replace the smoother based on Jacobi iteration with other iteration methods
(e.g., Gauss--Seidel iteration).
When applying this standard relaxation strategy, we found that replacing
the residual calculation (\ref{residual-L}) with $r^h = {f}^h - A^h
(E^hu^h)$ helps to maintain the desired convergence rate of the
V-Cycle scheme.
This extra extension has no effect in the continuous setting but seems
to improve the convergence rate of the discrete algorithm.

Even with this modified residual calculation, we sometimes observed
slightly faster convergence results using the Ruuth--Merriman smoothing
strategy. Thus for the remainder of this work, we use the Ruuth--Merriman smoothing strategy.

\subsection{Numerical examples}
\label{numerical-tests-multigrid}
We test our V-Cycle solver for the shifted Poisson's equation
on several closed curves in $\mathbb{R}^2$ and surfaces in $\mathbb{R}^3$.
We use cubic (degree 3) interpolations to construct $E^h$'s on each grid level.
In each round of the V-Cycle, $\nu_1 = 3$ pre-smoothings and $\nu_2 = 3$ post-smoothings are applied.
For the restriction and prolongation operators, linear (degree 1) interpolations are used
for the closest point extensions.
The linear system on the coarsest grid level is solved by Matlab's backslash.
The stopping criteria for all the tests are $\|x^{k+1}-x^k\|_{\infty}/\|x^k\|_{\infty}<10^{-6}$,
where $x^k$ and $x^{k+1}$ are the iterative solutions after the $k$-th and $(k+1)$-th rounds of V-Cycles.

\subsubsection{Unit circle}
We consider the simplest example from Section~\ref{poisson-circle}, where we solve
the shifted Poisson's equation $-\Delta_{\mathcal{S}}u + u = f$ on the unit circle.
The exact solution is chosen as $\sin(\theta)+\sin(12\theta)$.

\begin{figure}[htbp]
\begin{centering}
\subfigure[Decrease of relative residuals]{
\includegraphics[width=0.47\textwidth]{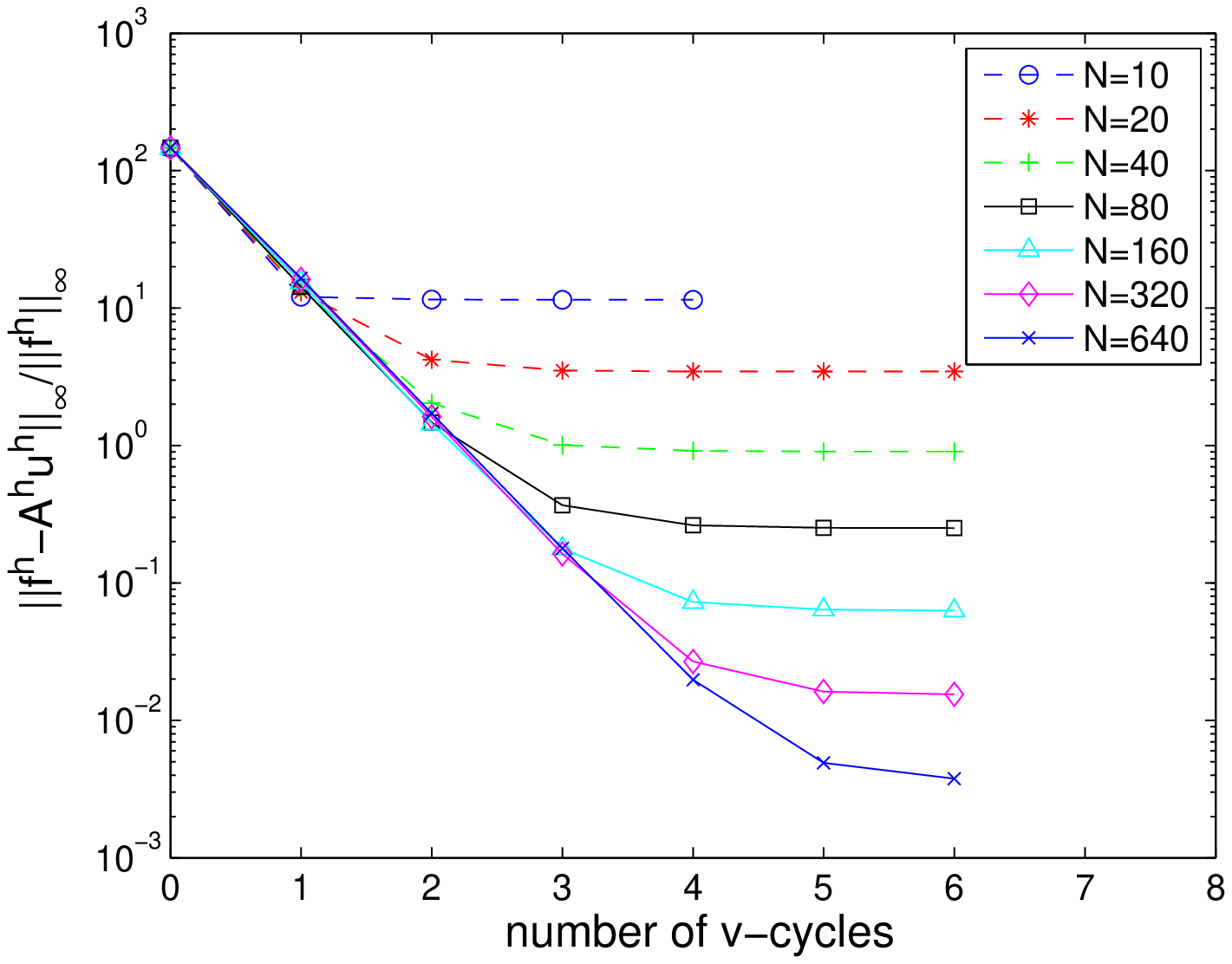}
}
\subfigure[Decrease of relative errors]{
\includegraphics[width=0.47\textwidth]{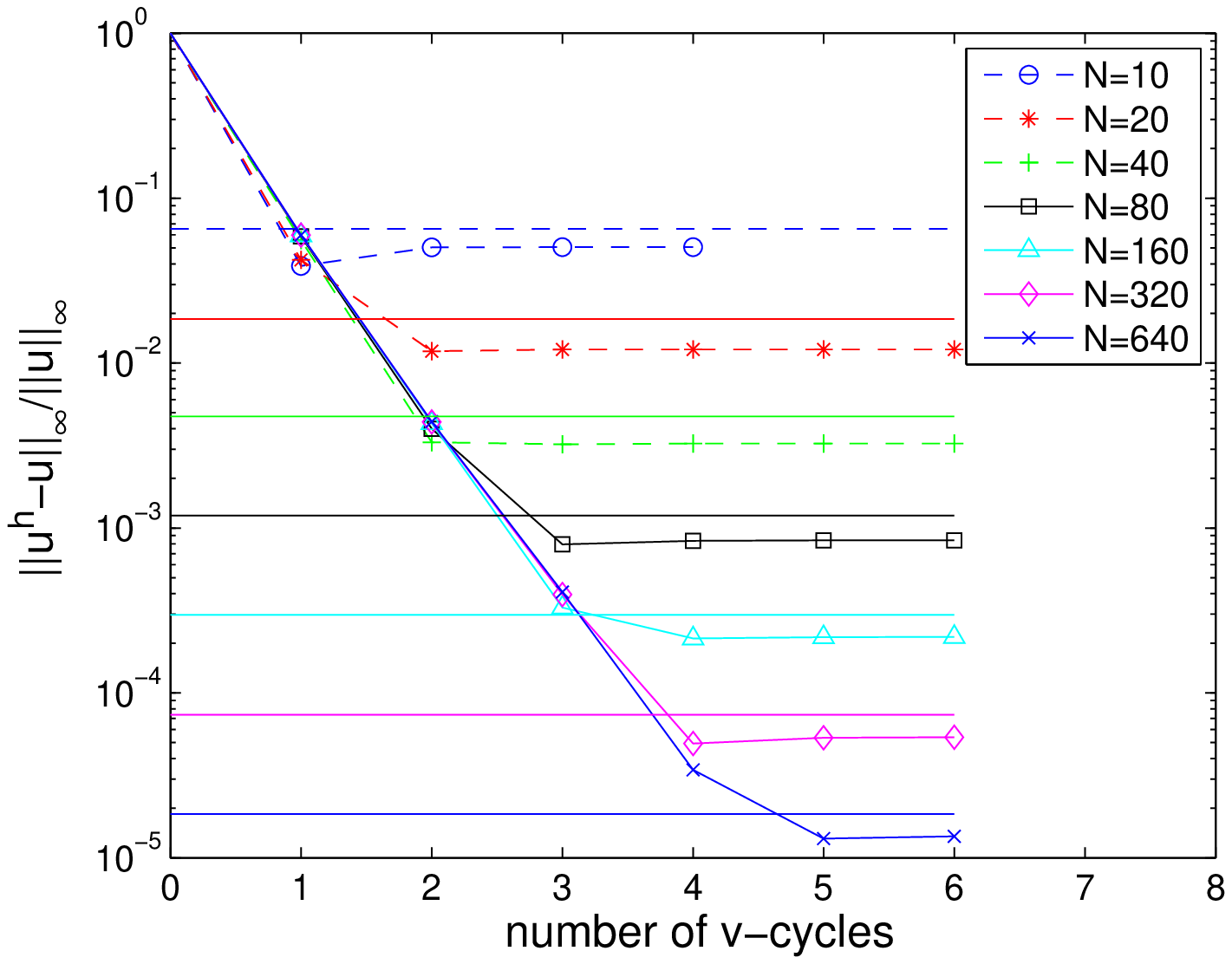}
}
\caption[Relaxation by $\tilde{A}^h=cI^h-L^h$ for the shifted Poisson's equation on the unit circle]{\label{circleL}
Convergence results on a circle, where the curves show a sequence of progressively finer grids, each solved with the V-Cycle multigrid method.
Each curve shows how the relative residuals/errors decrease as the number of V-Cycles increases.
Here $N=\frac{1}{\Delta x}$, where $\Delta x$ is the grid size of the finest mesh for that curve.
The coarsest grid in each case corresponds to $N=5$.
The horizontal lines in (b) are the numerical errors given by
Matlab's backslash. }
\end{centering}
\end{figure}

To test the convergence rate of our multigrid solver, we keep the coarsest grid size unchanged,
and reduce the finest grid size.
Figure~\ref{circleL}(a) shows the decreasing rate of the residuals of the algebraic linear systems,
and Figure~\ref{circleL}(b) shows the decreasing rate of the numerical errors.
Before the numerical solutions achieve the stopping criterion, they have already reached
the limit of the discretization errors, and the errors will not decrease after that;
so in Figure~\ref{circleL}(b) as the number of V-Cycles increases,
the errors first decrease and then stay almost the same.
This same behaviour is also observed in the residuals in Figure~\ref{circleL}(a):
as the number of V-Cycles increases, the residuals initially decrease.  However, because we
introduce additional discretization errors to the original linear system $A^hu^h=f^h$ (by performing discrete closest point extensions after each smoothing step),
the residuals eventually stagnate at the level of this additional discretization error.
In both Figures~\ref{circleL}(a) and~\ref{circleL}(b), the rate of decrease (the slope of the lines) stays almost the same as we refine the finest mesh, which is expected and desirable in a multigrid algorithm.

\subsubsection{Bean-shaped curve}
\begin{figure}[htbp]
\begin{center}
\subfigure[Numerical solution]{
	\label{bean-multigrid:a}
	\includegraphics[width=0.47\textwidth]{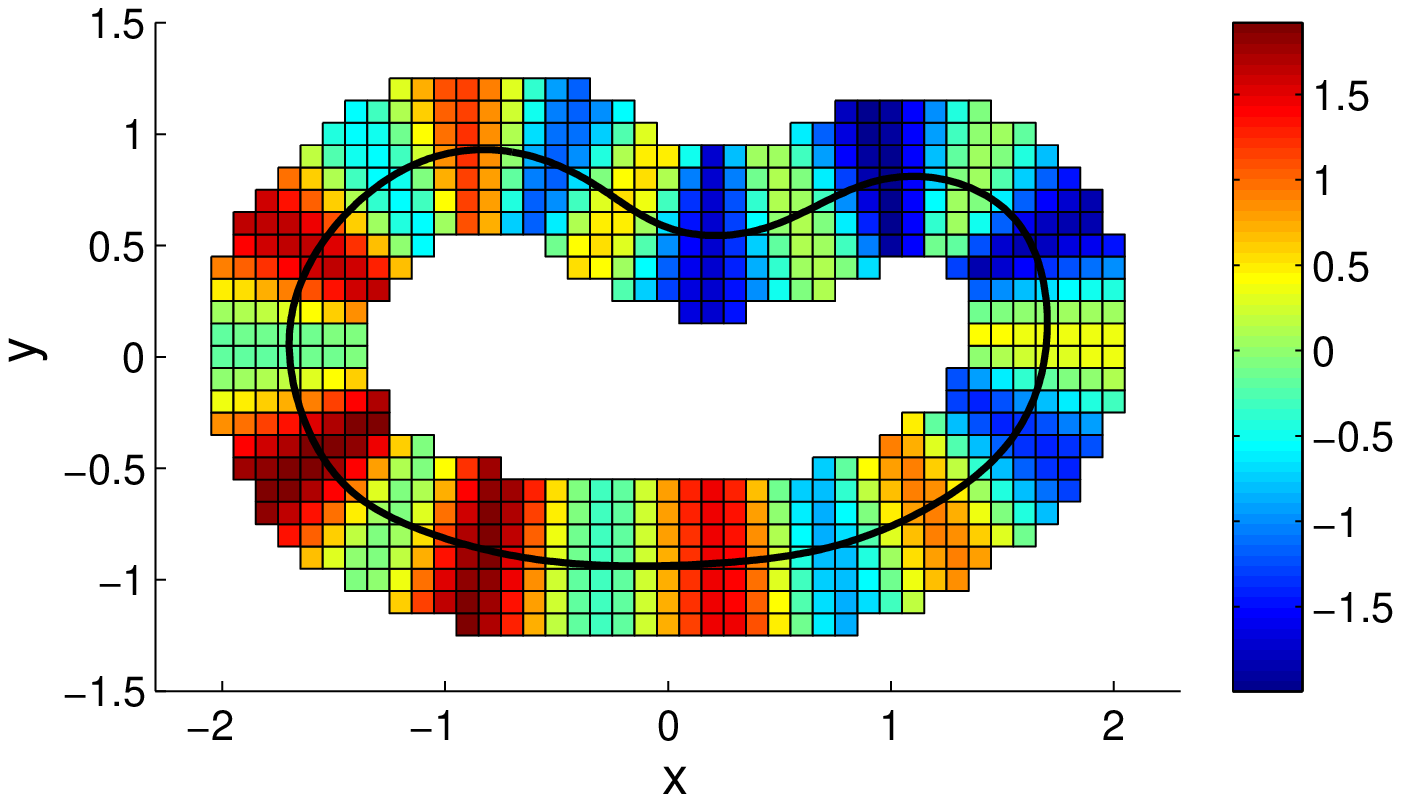}}
\subfigure[Decrease of relative errors]{
	\label{bean-multigrid:b}
	\includegraphics[width=0.47\textwidth]{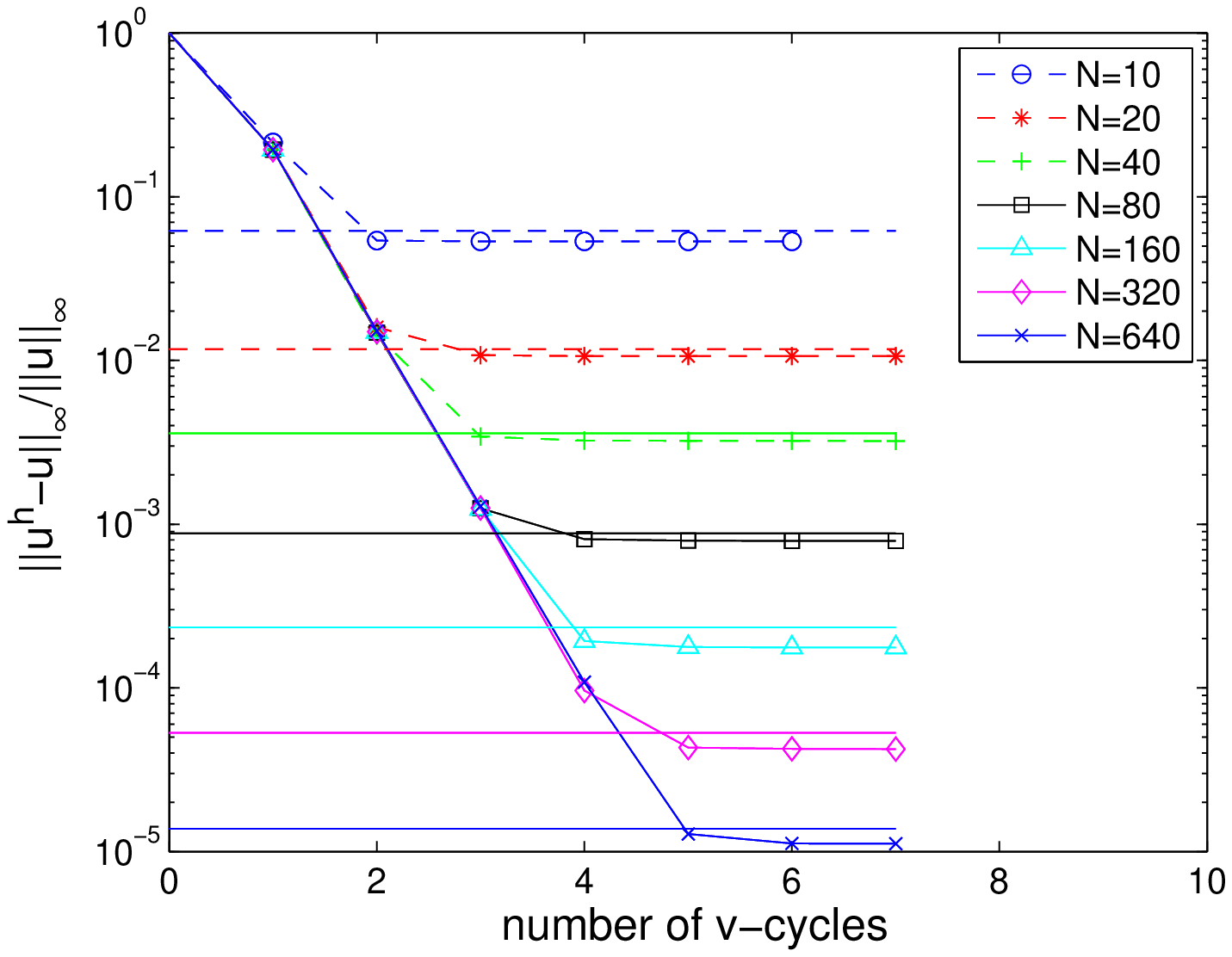}
	}
\caption{\label{bean-multigrid}
V-Cycle convergence results on a bean shaped curve. The horizontal lines in figure (b) are the numerical errors given by Matlab's backslash.}
\end{center}
\end{figure}

We test the shifted Poisson's equation on the bean-shaped curve as in Section~\ref{bean-backslash},
and the same exact solution is chosen for convergence studies.
Figure~\ref{bean-multigrid}(a) shows the numerical solution for the embedding equation with meshsize $\Delta x=0.1$.
Figure~\ref{bean-multigrid}(b) shows the convergence results of the multigrid method.
Note again that the convergence rate (slope of lines) is essentially independent of the finest grid size.

\subsubsection{Unit sphere}
\label{sphere}

To make an exact solution and perform a convergence study, we parameterize a sphere in the spherical coordinates $(\theta, \phi, r)$:
$x = r\sin\phi\cos\theta$, $y = r\sin\phi\sin\theta$, $z = r\cos\phi$,
where $\theta$ is the azimuth ranging from 0 to $2\pi$, $\phi$ is the
elevation from 0 to $\pi$, and $r$ is the radius.

\begin{figure}[htbp]
\begin{center}
\subfigure[Numerical solution]{
	\label{hemisphere-neumann-multigrid:a}
	\includegraphics[width=0.47\textwidth]{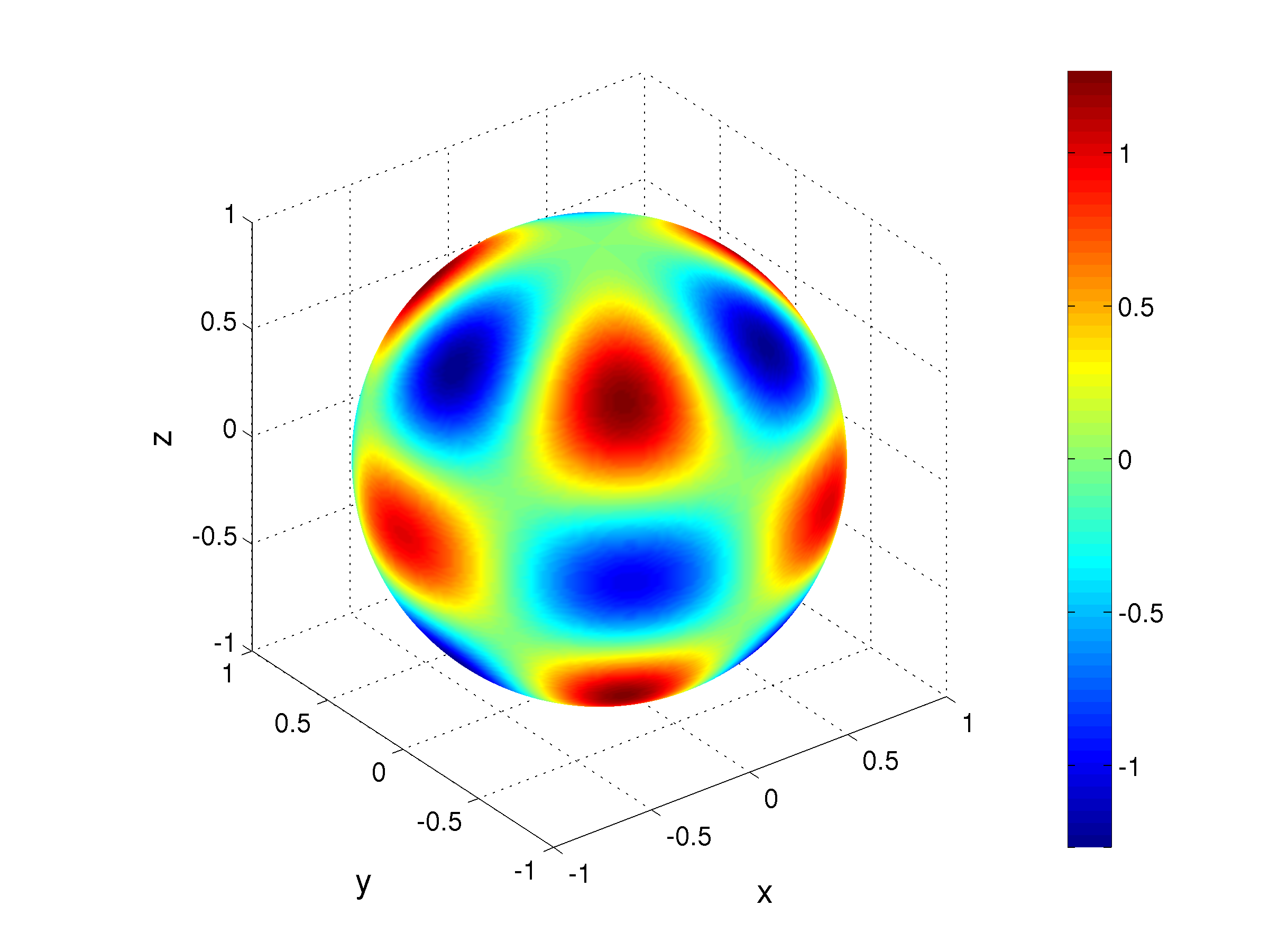}}
\subfigure[Decrease of relative errors]{
	\label{hemisphere-neumann-multigrid:b}
	\includegraphics[width=0.47\textwidth]{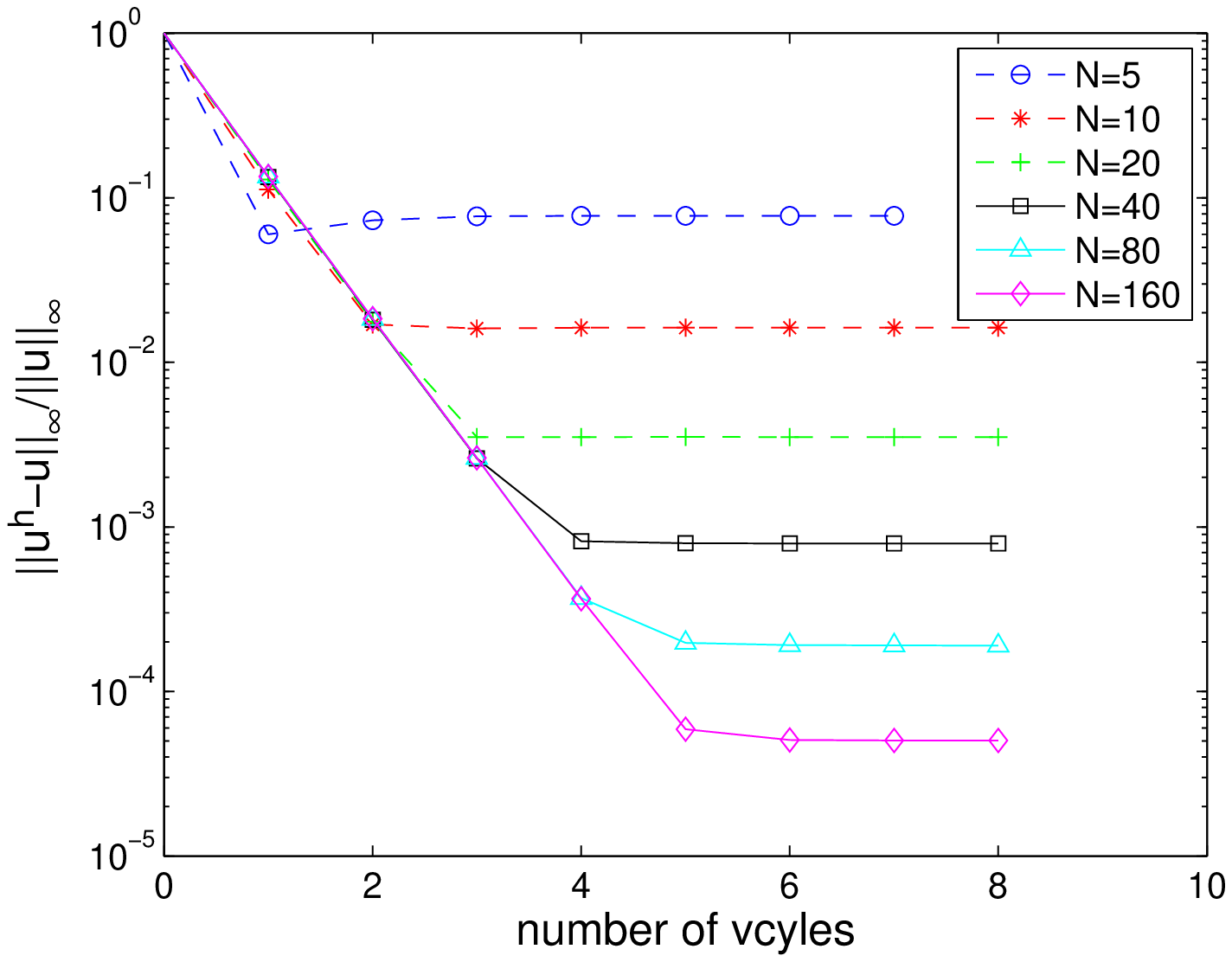}
	}
\caption{\label{hemisphere-neumann-multigrid}
V-Cycle convergence results on a unit sphere.}
\end{center}
\end{figure}

We emphasize this parametrization is not used in our algorithm: it is only used to compute the error.
For simplicity, we solve the equation on a unit sphere centered at the origin, and choose the exact solution as
the spherical harmonic $u(\theta,\phi)=\cos(3\theta)\sin(\phi)^3(9\cos(\phi)^2-1)$. The right-hand side function is
 $-\Delta_{\mathcal{S}}u+u=-29u(\theta,\phi)$. Figure~\ref{hemisphere-neumann-multigrid}(a)
shows the numerical solution for the embedding equation with grid size $\Delta x=0.1$,
visualized using the parametrization.
Figure~\ref{hemisphere-neumann-multigrid}(b) shows the convergence results of the multigrid method.

Besides convergence studies, we also investigate the computational cost including the degrees of freedom and matrix storage
for the discrete linear system (see Table~\ref{nnz}). Because our 3D computation is performed in a narrow band with
bandwidth proportional to the mesh size $\Delta x$, the number of degrees of freedom scales like $O(\Delta x^{-2})$, just as
if we performed a finite difference method to solve the Poisson's equation in $\mathbb{R}^2$. Since we use the standard 7-point stencil,
the number of nonzero entries (nnz) of the Laplacian matrix $L^h$ is about 7 times the number of degrees of freedom. Since we use degree
$3$ polynomial interpolation in a dimension by dimension fashion, the number of nonzero entries of the interpolation matrix
$E_3^h$ is about $4^3=64$ times the number of degrees of freedom.
The matrix $\tilde{A}^h$ used in our Ruuth--Merriman smoothing
approach is simply the diagonally shifted Laplacian so
$\text{nnz}(\tilde{A}^h) = \text{nnz}(L^h)$.
From the analysis in Section~\ref{convergence-proof},
$M^h$ has nearly the same sparsity pattern as $E_3^h$, except that $M^h$ might have more diagonal entries than $E_3^h$,
so the number of nonzero entries of $M^h$ is about $64$ or $65$ times the number of degrees of freedom.

\begin{table}[htbp]
\caption[Computational cost for the Closest Point Method on a sphere]{
\label{nnz} Computational cost for the Closest Point Method on a sphere. As
we decrease the mesh size by a factor of 2, the number of degrees of freedom (DOFs, i.e., $\text{length}(u^h)$)
increases by a factor of 4. $\text{nnz}(L^h)\approx 7\times\text{DOFs}$,
$\text{nnz}(E^h)\approx 64\times\text{DOFs}$, $\text{nnz}(M^h)\approx 64.5\times\text{DOFs}$.
}
\centering
\begin{small}
\begin{tabular}{ccccc}
\hline
\rule[-0.9ex]{0pt}{3.2ex}  
$\Delta x$ & $\text{length}(u^h)$ & $\text{nnz}(L^h)$ & $\text{nnz}(E_3^h)$ & $\text{nnz}(M^h)$ \\
\hline
0.2        & 3190        & 20758            & 204160           & 205686   \\
0.1        & 10906       & 71962            & 697984           & 702714   \\
0.05       & 41870       & 277358           & 2679680          & 2697038  \\
0.025      & 166390      & 1103734          & 10648960         & 10717230  \\
0.0125     & 663454      & 4402366          & 42461056         & 42731646 \\
0.00625    & 2651254     & 17593318         & 169680256        & 170760222 \\
\hline
\end{tabular}
\end{small}
\end{table}

We also compare the CPU times for solving the linear system by Matlab's backslash and
our V-Cycle multigrid scheme.
Table~\ref{cpu-sphere} shows the performance comparison.
Computations were performed on a 3.47 GHz Intel Xeon X5690 6-core processor,
with 100 GB memory.
For the smallest problem, Matlab's backslash works
 better than multigrid; but for larger problems, multigrid works better. What is more, the CPU time
 of multigrid scales roughly as $O(\Delta x^{-2})$, i.e., proportional to the number of degrees of freedom, which
 is optimal and what we expect for a multigrid scheme \cite{plain:briggs}.

\begin{table}[htbp]
\begin{center}
\begin{minipage}{6cm}
\begin{small}
\begin{tabular}{ccc}
\hline
$\Delta x$  & Backslash        & V-Cycle  \\
\hline
0.2         & 0.169           & 0.294 \\
0.1         & 1.046           & 0.551 \\
0.05        & 6.967           & 1.326 \\
0.025       & 59.94           & 4.377 \\
0.0125      & 466.5           & 16.66 \\
0.00625     & N/A            & 67.42 \\
\hline
\end{tabular}
\end{small}
\end{minipage}
\begin{minipage}{6cm}
\caption[CPU time for a sphere]{
\label{cpu-sphere}
CPU time for solving the linear system arising from the
Poisson's equation on a unit sphere with the Closest Point Method.
Matlab's backslash versus V-Cycle multigrid method.
}
\end{minipage}
\end{center}
\end{table}

\subsubsection{Torus}
Again, to calculate an exact solution, we begin with a parametrization of a torus:
$x(\theta,\phi) = (R+r\cos \phi)\cos \theta$,
$y(\theta,\phi) = (R+r\cos \phi)\sin \theta$,
$z(\theta,\phi) = r\sin\phi$,
where $\phi$, $\theta$ are in the interval $[0,2\pi)$, $R$ is the distance from
the center of the tube to the center of the torus, and $r$ is the radius
of the tube.

\begin{figure}[htbp]
\begin{center}
\subfigure[Numerical solution]{
	\label{torus-multigrid:a}
	\includegraphics[width=0.47\textwidth]{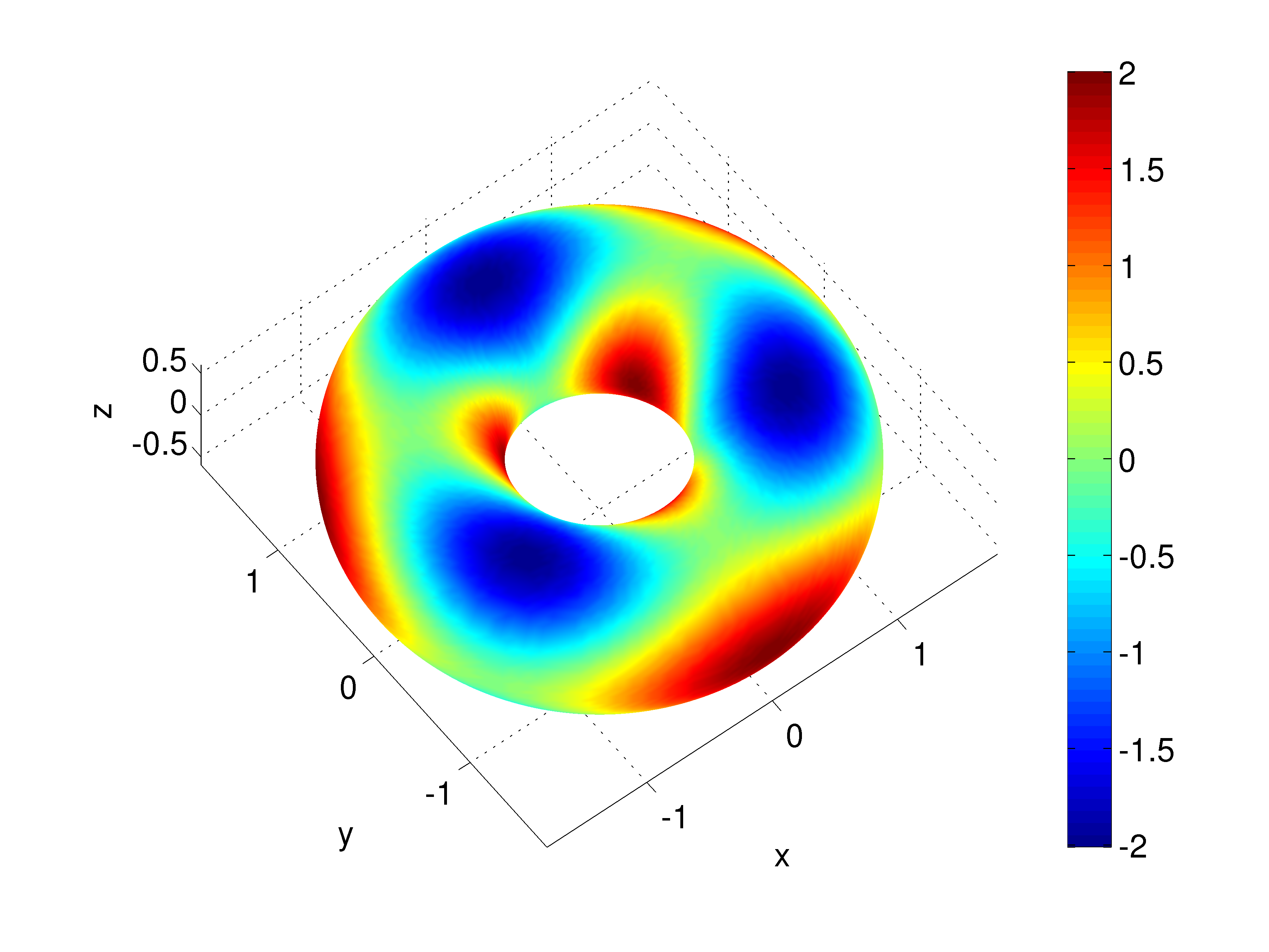}}
\subfigure[Decrease of relative errors]{
	\label{torus-multigrid:b}
	\includegraphics[width=0.47\textwidth]{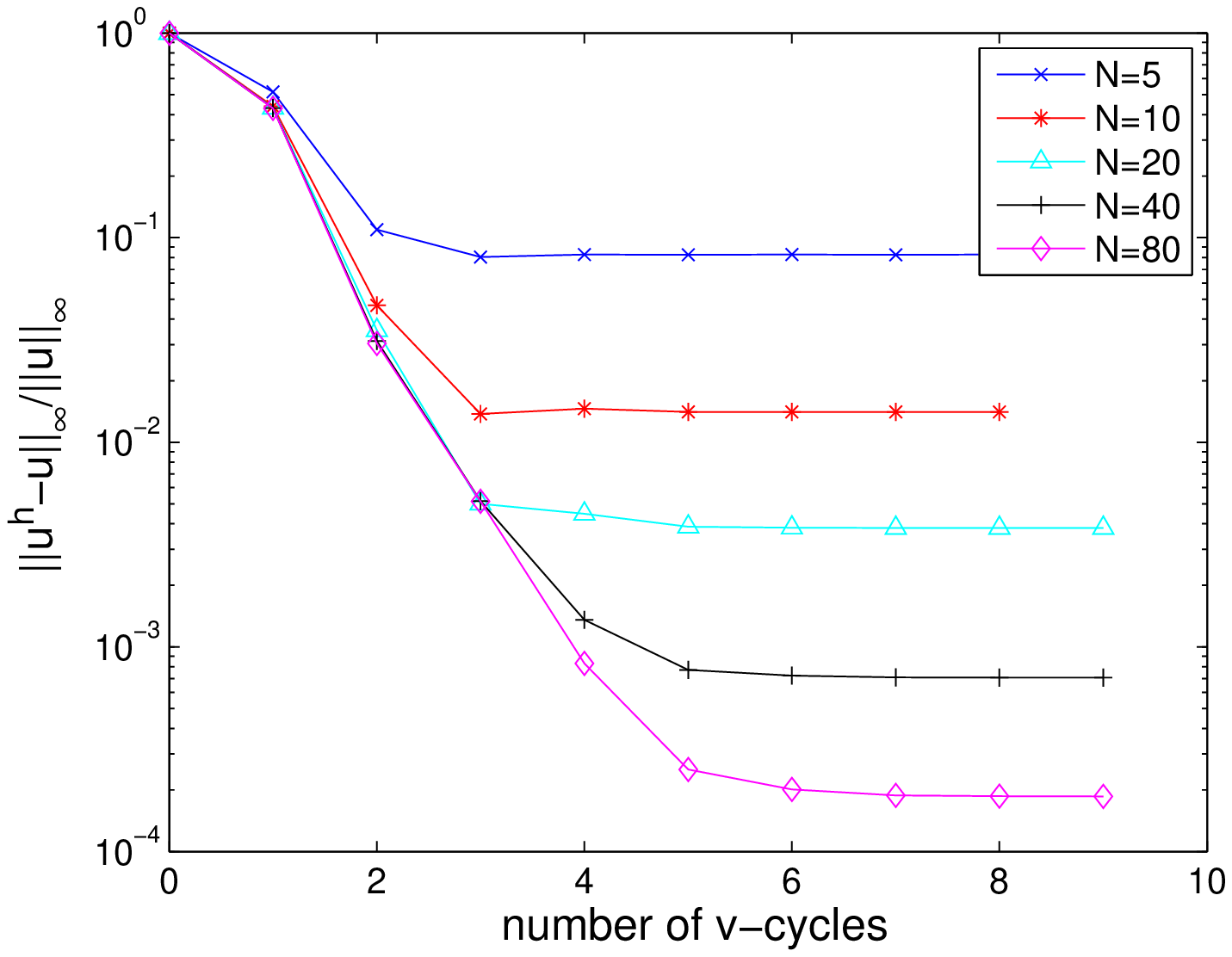}
	}
\caption{\label{torus-multigrid}
V-Cycle convergence results on a torus.}
\end{center}
\end{figure}

We test the shifted Poisson's equation $-\Delta_{\mathcal{S}} u + u=f$ on a torus with $R=1.2$, and $r=0.6$.
The exact solution is chosen as $u(\theta,\phi) =\sin(3\theta)+\cos(2\phi)$, and
the right hand side function is
$f(\theta,\phi) = \frac{9\sin(3\theta)}{(R+r\cos\phi)^2} - \frac{2\sin\phi\sin (2\phi)}{r(R+r\cos\phi)}
+\frac{4\cos(2\phi)}{r^2} + u(\theta,\phi)$.
Figure~\ref{torus-multigrid} shows the numerical solution and the multigrid convergence results.

\subsubsection{A level set example} \label{level-set-test}
\begin{figure}[htbp]
\begin{center}
\subfigure[Numerical solution]{
	\label{ls-mg:a}
	\includegraphics[width=0.47\textwidth]{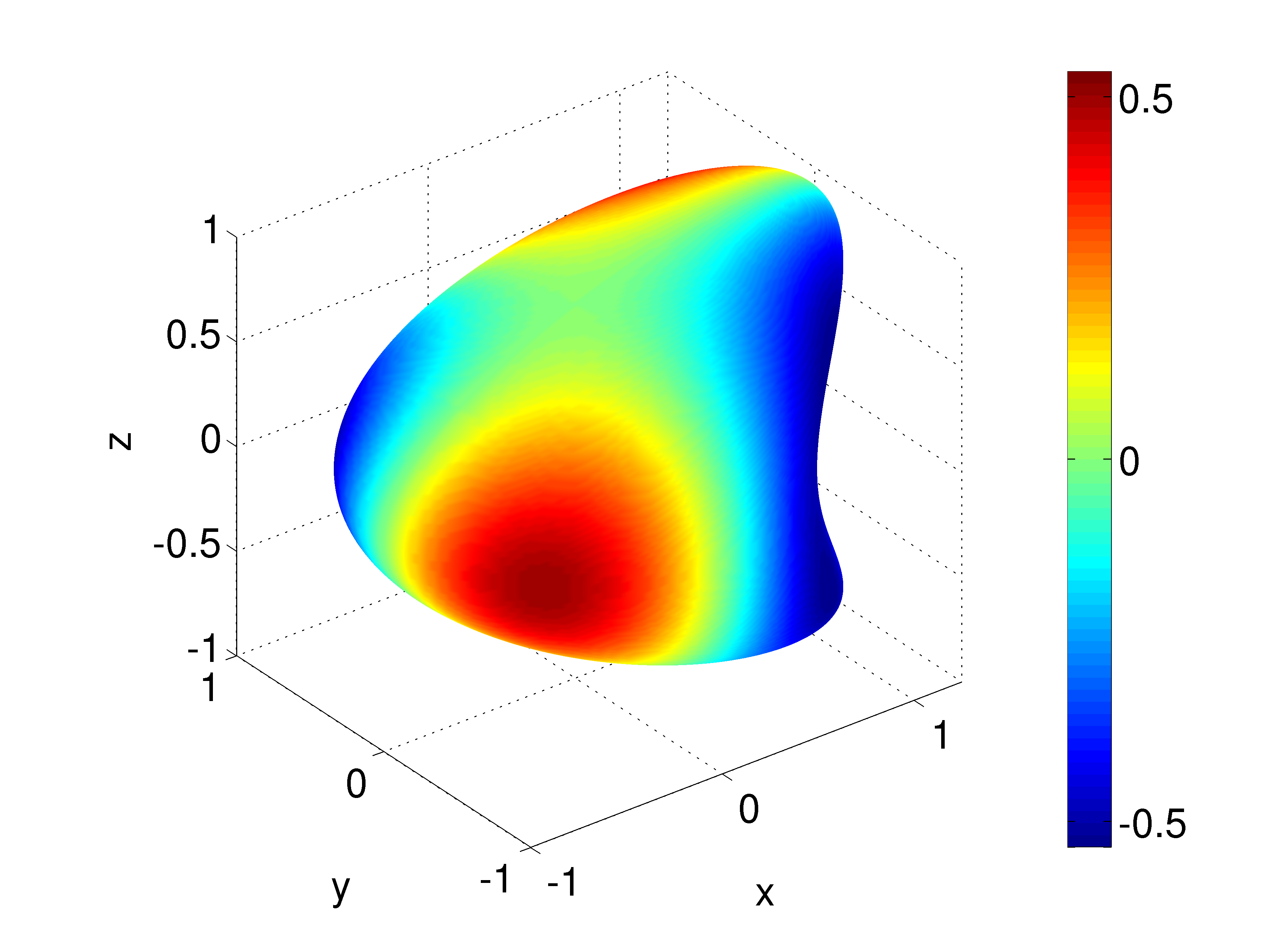}}
\subfigure[Decrease of relative errors]{
	\label{ls-mg:b}
	\includegraphics[width=0.47\textwidth]{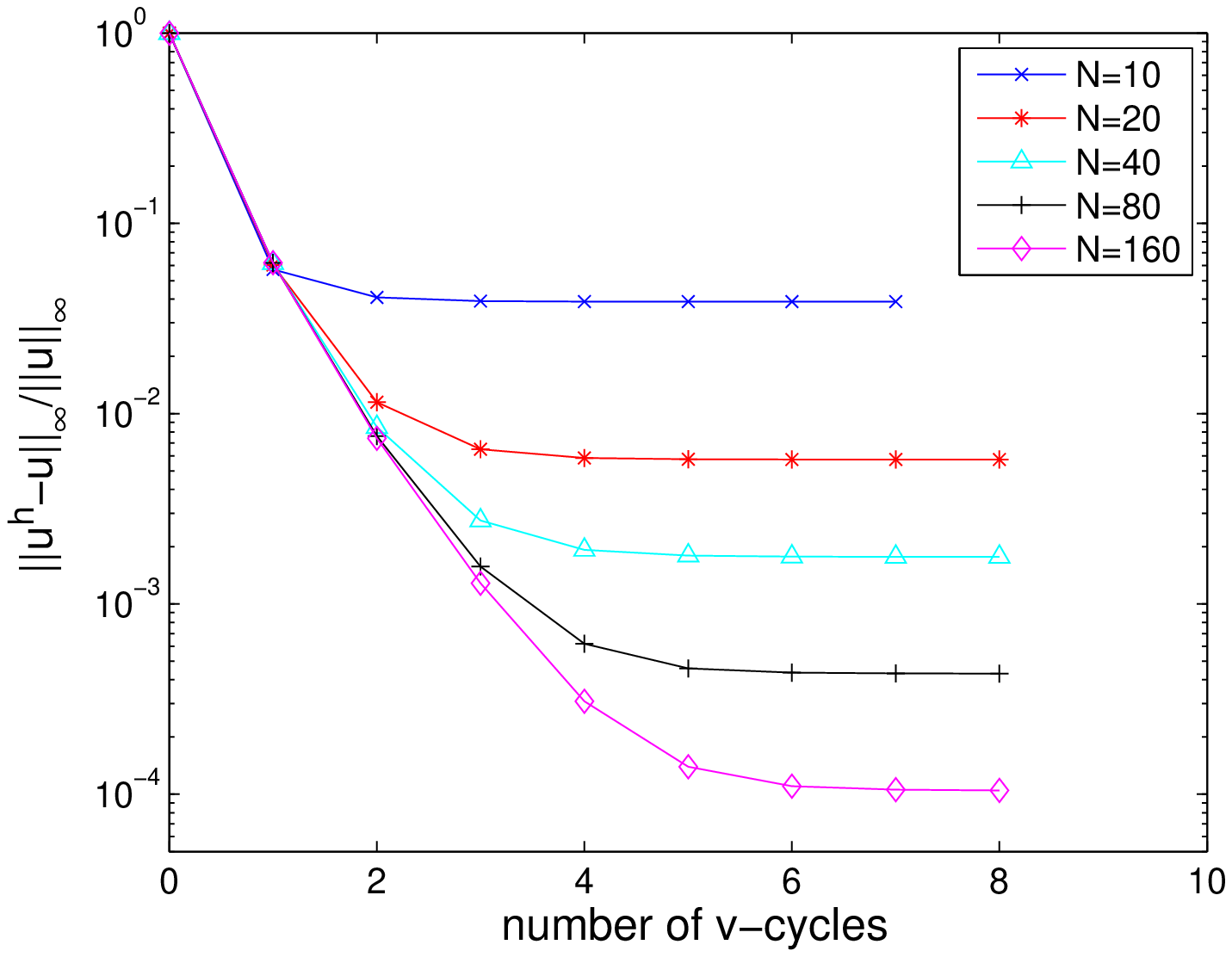}
	}
\caption{\label{ls-mg}
V-Cycle convergence results on the Dziuk surface.}
\end{center}
\end{figure}

Following Dziuk \cite{plain:sfem-error}, we consider the surface defined by
$\mathcal{S} = \{ \bm{x}\in \mathbb{R}^3 ~|~ (x_1-x_3^2)^2+x_2^2+x_3^2=1\}$,
and the shifted Poisson's equation $-\Delta_{\mathcal{S}}u+u = f$ on $\mathcal{S}$.
The exact solution is chosen to be $u(\bm{x})=x_1x_2$. We compute the right-hand side
function by $f=-\nabla_{\mathcal{S}}\cdot \bm{v}$, where $\quad \bm{v} = \nabla_{\mathcal{S}}u = \nabla u - (\nabla u\cdot \bm{n})\bm{n}$,
and $\nabla_{\mathcal{S}}\cdot \bm{v} = \nabla\cdot \bm{v} - \sum_{j=1}^3(\nabla \bm{v}_j\cdot \bm{n})\bm{n}_j$.
Here $\bm{n}$ is the normal vector, and in this case,
$\bm{n}(\bm{x}) = \left(~ x_1-x_3^2, ~ x_2, ~ x_3(1-2(x_1-x_3^2)) ~ ) ~/~ (1+4x_3^2(1-x_1-x_2^2))^{1/2} ~\right)$.
For each grid point, we compute the closest point on $\mathcal{S}$
by minimizing the squared distance function using Newton's method.
Figure~\ref{ls-mg} shows the numerical solution and multigrid convergence results.

\subsubsection{Elliptic equation with variable diffusion coefficients on a sphere}
\label{section-sphere-vc}
Finally, we consider an example with scalar variable diffusion coefficients,
$-\nabla_{\mathcal{S}}\cdot(a(\bm{y})\nabla_{\mathcal{S}}u(\bm{y}))+cu(\bm{y})=f(\bm{y})$,
for $\bm{y}\in\mathcal{S}$.

To formulate the continuous embedding equation, we simply choose $\bm{A}(\text{cp}(\bm{x}))=a(\text{cp}(\bm{x}))$
in the general embedding equation (\ref{embedding-general}), and get
$-[\nabla\cdot(a(\text{cp}(\bm{x}))\nabla\tilde{u})](\text{cp}(\bm{x})) + c\tilde{u}(\bm{x}) + \gamma(\tilde{u}(\bm{x})-\tilde{u}(\text{cp}(\bm{x})))
= \tilde{f}(\bm{x})$
for $\bm{x}\in B(\mathcal{S})$.
Similar to the Poisson case before, this embedding equation
can be discretized using standard finite differences and interpolation schemes, see also \cite{plain:mol}.
The matrix formulation of the discretization is
$-E^h \tilde{L}^h u^h + \gamma(I^h-E^h)u^h = f^h$,
where
$\tilde{L}^h =  \bm{D}_x^b \text{diag}(\bm{A}_x^f a^h) \bm{D}_x^f
         + \bm{D}_y^b \text{diag}(\bm{A}_y^f a^h) \bm{D}_y^f
         + \bm{D}_z^b \text{diag}(\bm{A}_z^f a^h) \bm{D}_z^f$.
Here we have used similar notations to those in \cite{plain:mol}.
$\bm{D}_x^b$ is the matrix corresponding to backward differences in the $x$ direction,
$\bm{D}_x^f$ is the matrix corresponding to forward differences in the $x$ direction;
$\bm{D}_y^b$, $\bm{D}_y^f$, $\bm{D}_z^b$ and $\bm{D}_z^f$ have similar meanings. The column vector $a^h$ has
the $i$-th entry equal to $a(\text{cp}(\bm{x}_i))$, where $\bm{x}_i$ is the $i$-th
grid point. $\bm{A}_x^f$ is the forward averaging matrix corresponding to calculating the average of
two neighboring values of $a^h$ along the positive $x$ direction. $\bm{A}_x^f a^h$ is a column vector
and $\text{diag}(\bm{A}_x^f a^h)$ is a diagonal matrix with diagonal entries equal to $\bm{A}_x^f a^h$;
$\text{diag} (\bm{A}_y^f a^h)$ and $\text{diag} (\bm{A}_z^f a^h)$ have similar meanings.

In the multigrid algorithm, we continue to use the Ruuth--Merriman style relaxation strategy,
except changing $L^h$ to $\tilde{L}^h$.

\begin{figure}[htbp]
\begin{center}
\subfigure[Decrease of relative residuals]{
	\includegraphics[width=0.47\textwidth]{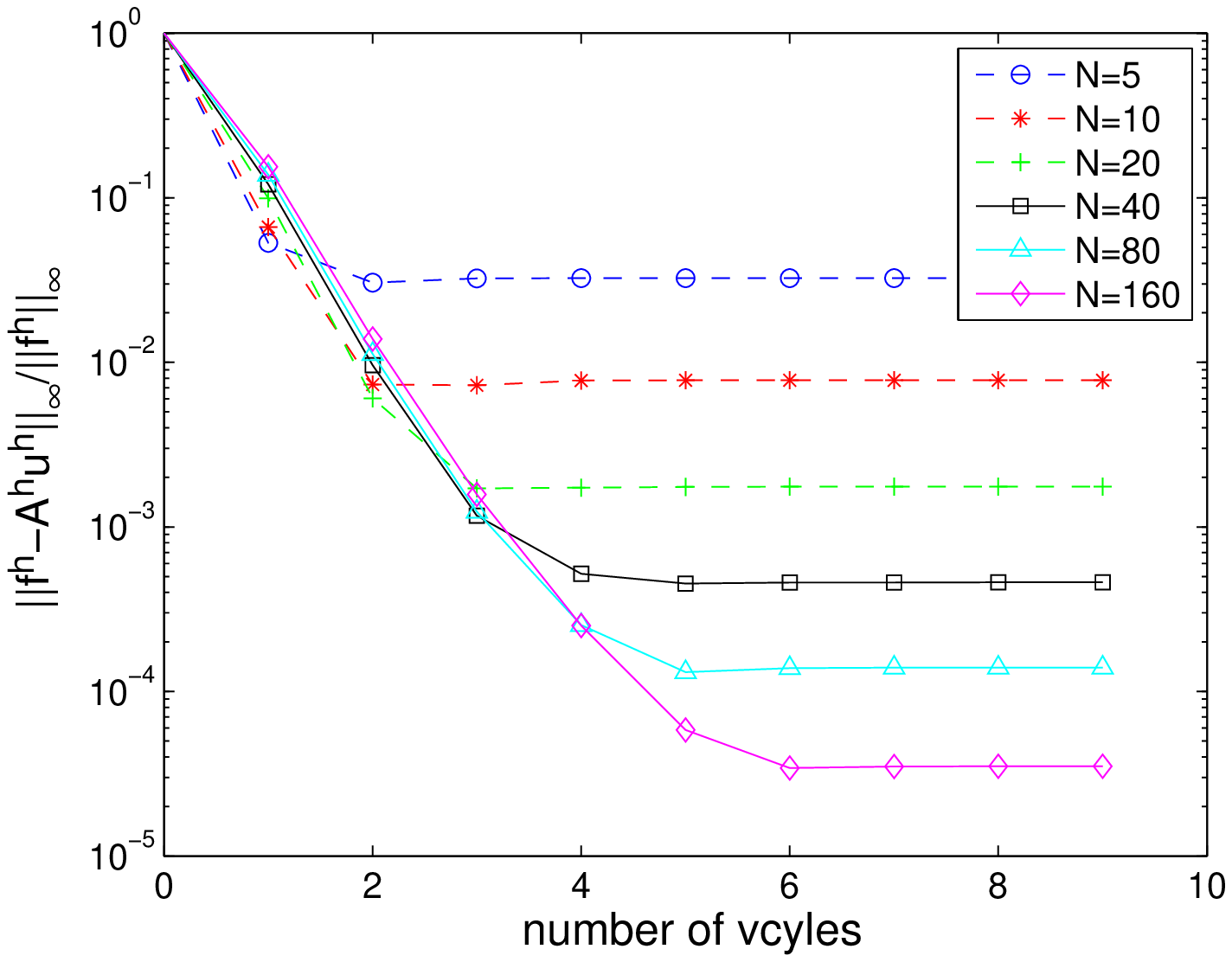}}
\subfigure[Decrease of relative errors]{
	\includegraphics[width=0.47\textwidth]{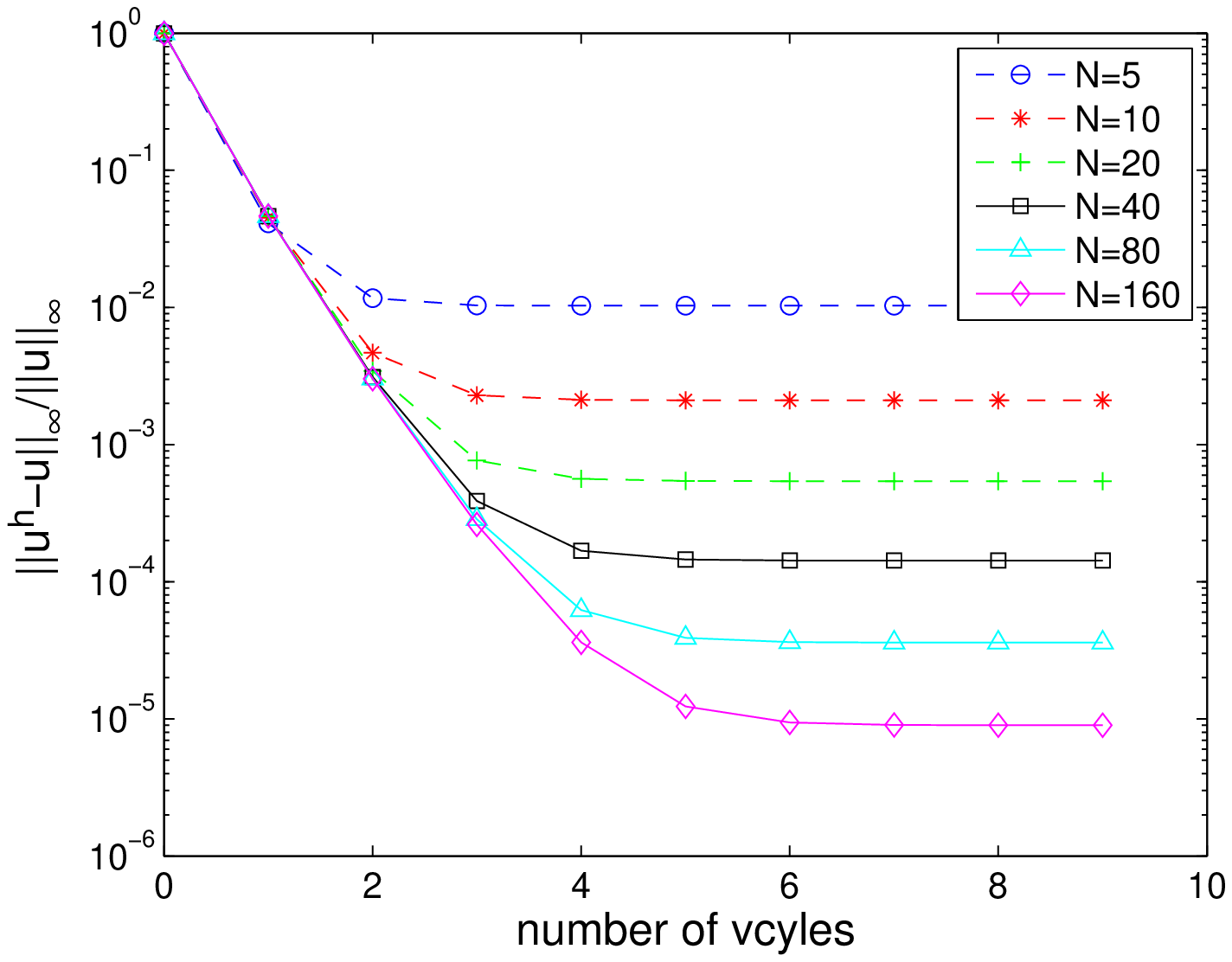}
	}
\caption{\label{sphere-vc-figure}
V-Cycle convergence results for a variable diffusion coefficient example.}
\end{center}
\end{figure}

We test on a unit sphere. We choose the exact solution
$u(\phi)=\cos\phi$, diffusion coefficient $a(\phi)=\cos\phi+1.5$,
and the right hand side function is $f(\phi) = 2\cos\phi(\cos\phi+1.5)-(\sin\phi)^2+u(\phi)$.
The multigrid convergence results are shown in Figure~\ref{sphere-vc-figure}
and look similar to the Poisson problems shown earlier.

\section{Conclusions}
\label{section-conclusion}

We adapt the Closest Point Method to solve elliptic PDEs on surfaces.
In particular, we analyze the shifted Poisson's equation $-\Delta_{\mathcal{S}}u+cu=f$ in detail.
We formulate an embedding equation which includes a \emph{consistency condition} from the surface PDE
and a \emph{side condition} to enforce the solution to be constant along the normals to the surface.
We then discretize the embedding equation using standard centered finite difference methods
and Lagrange interpolation schemes. We prove the convergence for the difference scheme
in a simple case where the underlying surface is a closed curve embedded in $\mathbb{R}^2$.

We then propose a specific geometric multigrid method to solve the resulting
large sparse linear system. The method makes full use of the
closest point representation of the surface and uniform Cartesian
grids in different levels.
Numerical tests suggest that the convergence speeds do
not deteriorate as we increase the number of grid levels while fixing the size of the coarsest grid,
just as if we performed a standard multigrid algorithm for the (Cartesian) Laplacian operator.

There are many other possible multigrid methods; in particular,
Section~\ref{section-standard-relax} mentions an alternative approach based directly on our modified equation (which additively combines the side and consistency conditions).
Future work could investigate this further.


Numerical results in Section~\ref{section-sphere-vc} show that our approach also works for non-constant coefficient elliptic PDEs
with scalar diffusivity $a(\bm{y})$. 
It would be useful for some applications to investigate more general elliptic PDEs,
replacing $a(\bm{y})$ with a symmetric and semi-positive definite diffusivity matrix $\bm{A}(\bm{y})$.
Our approach should work here too although we have not tested it.

Finally, our analysis of the spatial discretizations of elliptic operators
can be used to analyze embedded method-of-lines
approaches for parabolic PDEs on general surfaces \cite{plain:mol}.

\bibliographystyle{siam}
\bibliography{mybib}

\end{document}